\documentclass{amsart}
\usepackage{mathrsfs}
\usepackage{amsfonts}
\usepackage[all]{xy}
\usepackage{latexsym,amsfonts,amssymb,amsmath,amsthm,amscd,amsmath,graphicx}
\usepackage{tikz}
\newtheorem{thm}{Theorem}[section]
\newtheorem{cor}[thm]{Corollary}
\newtheorem{lem}[thm]{Lemma}
\newtheorem{prop}[thm]{Proposition}
\theoremstyle{definition}
\newtheorem{defn}[thm]{Definition}
\newtheorem{exam}[thm]{Example}
\theoremstyle{remark}
\newtheorem{rem}[thm]{Remark}
\numberwithin{equation}{section}

\newcommand{\p}{\partial}

\newcommand{\sF}{{\mathcal F}}
\newcommand{\sG}{{\mathcal G}}

\newcommand{\sK}{{\mathcal K}}
\newcommand{\sL}{{\mathcal L}}

\newcommand{\sO}{{\mathcal O}}

\newcommand{\sQ}{{\mathcal Q}}

\newcommand{\sW}{{\mathcal W}}

\newcommand{\A}{{\mathbb A}}

\newcommand{\C}{{\mathbb C}}

\newcommand{\F}{{\mathbb F}}
\newcommand{\G}{{\mathbb G}}

\newcommand{\N}{{\mathbb N}}

\newcommand{\Z}{{\mathbb Z}}

\newcommand{\x}{\xrightarrow}
\DeclareMathOperator{\DR}{DR}

\numberwithin{equation}{section}

\begin {document}
\topmargin= -.2in \baselineskip=20pt

\title{Composition series for GKZ-systems}
\author{Jiangxue Fang}
\address{Department of Mathematics,
Capital Normal University, Beijing 100148, P.R. China} \email{fangjiangxue@gmail.com}
\date{}

\maketitle

\begin{abstract}
In this paper, we find a composition series of GKZ-systems with semisimple successive quotients. We also study the composition series of the corresponding perverse sheaves and compare these two composition series under the Riemann-Hilbert correspondence.
\end{abstract}

\section{\textbf{Introduction}}

\subsection{GKZ-systems and their Fourier transforms} Let $A=(a_{ij})\in{\rm M}_{n\times N}(\Z)$ be a matrix with column vectors $a_1,\ldots,a_N$.
In the 1980s, Gelfand, Kapranov and Zelevinsky associated  to $A$ a class of $D$-modules which are now called  GKZ-systems or $A$-hypergeometric systems and are defined as follows.

Let $x_A=\{x_1,\ldots,x_N\}$ be a set of $N$ variables, and let $\C[x_A]=\C[x_1,\ldots,x_N]$. Let $\p_A=\{\p_1,\ldots,\p_N\}$ be the corresponding partial derivative operators on $\C[x_A]$, and let $\C[\p_A]=\C[\p_1,\ldots,\p_N]$. Let $\A^A={\rm Spec}\,\C[x_A]$, and let $D_A=\C[x_A,\p_A]$ be the Weyl algebra.

For any $v=(v_1,\ldots,v_N)^{\rm t}\in\Z^N$, define
$$\Box_v=\prod\limits_{v_j>0}\p_{j}^{v_j}-\prod\limits_{v_j<0}\p_{j}^{-v_j}.$$
For any parameter $\beta=(\beta_1,\ldots,\beta_n)^{\rm t}\in\C^n$, the $A$-hypergeometric system is the $D_A$-module
$$M_A(\beta)=D_A\Big/\sum\limits_{\substack{v\in\Z^N\\Av=0}}D_A\Box_v+\sum\limits_{i=1}^nD_A\Big(\sum\limits_{j=1}^Na_{ij}x_j\p_{j}-\beta_i\Big).$$

In this paper, all schemes are of finite type over $\C$ and most of them are affine. So, we don't distinguish between $D$-modules on an affine scheme and those over its coordinate ring.

The GKZ-system $M_A(\beta)$ is a holonomic $D_A$-module, by \cite[Theorem 3.9]{A} or Lemma \ref{lao}. A basic problem for $M_A(\beta)$ is how to find a canonical filtration with semisimple successive quotients. Such a filtration was first considered by Batyrev in \cite{B} to study the mixed Hodge structure of affine hypersurfaces in a torus, and it was shown by Stienstra in \cite{S} that the aforementioned filtration restricted to the generic fiber of the GKZ-system corresponds to the weight filtration of the relative cohomology group. Hence, it is natural to hope that this filtration extended to the whole space has semisimple successive quotients. In the slides of conference talk \cite{A2}, A. Adolphson defined a filtration on the GKZ-system $M_A(\beta)$ and stated that such a filtration has semisimple successive quotients. Unfortunately, the proof seems to have not been published up to now.

In this paper, we show that under certain conditions, the filtration on $M_A(\beta)$ considered in \cite{A2} has semisimple successive quotients (see Theorem \ref{baaa}). In Example \ref{hou}, we also show that such a filtration does not always have semisimple successive quotients without the aforementioned conditions. To do this, we only need to study such a filtration on the Fourier transform of $M_A(\beta)$ by the exactness of the Fourier transform (see Definition \ref{che}). We also study the filtration on the corresponding perverse sheaf of the Fourier transform of $M_A(\beta)$ under the Riemann-Hilbert correspondence (see Theorem \ref{yaode}).

In this paper, we emphasize that we do not assume that the rank of $A$ is $n$ or $\N A\cap-\N A=0$.
Let $\{t_1,\ldots,t_n\}$ be  a set of $n$ variables. For any $1\leq j\leq N$, set $t^{a_j}=\prod\limits\limits_{i=1}^nt_i^{a_{ij}}$. Let $S_A$ be the $\C$-subalgebra of $\C[t_1^{\pm1},\ldots,t_n^{\pm1}]$ generated by $t^{a_1},\ldots,t^{a_N}$.
The Fourier transform of $M_A(\beta)$ is
$$N_A(\beta)=D_A\Big/D_AI_A+\sum\limits_{i=1}^nD_A\Big(\sum\limits_{j=1}^Na_{ij}\p_{j}x_j+\beta_i\Big),$$
where $I_A$ is the kernel of the epimorphism
$$\C[x_A]\to S_A,\;x_j\mapsto t^{a_j}\hbox{ for any }1\leq j\leq N.$$
This epimorphism defines a closed immersion $ i_A\colon X_A={\rm Spec}\,S_A\to\A^A$.
We have an action $T_A\times X_A\to X_A$ of the torus $T_A={\rm Spec}\,\C[t^{\pm a_1},\ldots,t^{\pm a_N}]$ on $X_A$ defined by the homomorphism
$$\C[t^{a_1},\ldots,t^{a_N}]\to\C[t^{\pm a_1},\ldots,t^{\pm a_N}]\otimes_\C\C[t^{a_1},\ldots,t^{a_N}],\;\;t^{a_j}\mapsto t^{a_j}\otimes t^{a_j}.$$
In Lemma \ref{lao}, we will prove that $N_A(\beta)$ is a regular holonomic $D_A$-module supported on the toric variety $X_A$.

We have two methods to study filtrations on $N_A(\beta)$. The first one uses the Euler-Koszul complexes by combinatorial properties of the convex polyhedral cone in $\mathbb R^n$ generated by the column vectors of $A$. The second one uses the Riemann-Hilbert correspondence to study filtrations on the corresponding perverse sheaf. Finally, we compare these two filtrations via the de Rham functor.

\subsection{Functors on $D$-modules and perverse sheaves}
We introduce the following notation.

Let $X$ be a smooth scheme with a closed subscheme $Z$. We denote by ${\rm Mod}(D_X)$ the abelian category of left $D$-modules on $X$. The full subcategory of ${\rm Mod}(D_X)$ consisting of regular holonomic $D_X$-modules with support on $Z$ is denote by ${\rm Mod}_{rh}^Z(D_X)$.
Let $D^b(D_X)$ be the derived category of ${\rm Mod}(D_X)$.
The full subcategory of $D^b(D_X)$ consisting of objects whose cohomology are $\sO_X$-quasi coherent (resp. regular holonomic with support on $Z$) is denoted by $D^b_{qc}(D_X)$ (resp. $D^{b,\,Z}_{rh}(D_X)$).

Given a morphism $f\colon X\to Y$ of smooth schemes, denote by $D_{X\to Y}$ and $D_{Y\leftarrow X}$ the transfer bimodules. Given $M \in D^b(D_X)$ and $N\in D^b(D_Y )$, the direct and inverse images for $D$-modules are defined by
\begin{eqnarray*}
f_+M:=Rf_*(D_{Y\leftarrow X}\otimes^L_{ D_X} M)\hbox{ and }f^+N:=D_{X\to Y}\otimes^L_{ f^{-1}D_Y}f^{-1} N.
\end{eqnarray*}
Here, I replace the notation $Lf^*$ and $\int_f$ in \cite[Chapter 1, 1.5]{HTT} by $f^+$ and $f_+$, respectively. By \cite[Theorem 6.1.5]{HTT}, $f_+$ and $f^+$ preserve quasi-coherence and regular holonomicity.

If $f\colon X\to Y$ is an open immersion of smooth schemes, the minimal extension $f_{!+}M$ of a regular holonomic $D_X$-module $M$ is defined to be the smallest $D_Y$-submodule of $H^0f_+M$ whose restriction on $X$ coincides with $M$. More generally, if $f$ is an immersion (that is, $f=j\circ i$ for some closed immersion $i$ and open immersion $j$), define $f_{!+}M=j_{!+}(i_+M)$. The regular holonomic $D_Y$-module $f_{!+}M$ depends only on $f$ and $M$.
This definition of minimal extension coincides with \cite[Definition 3.4.1]{HTT}.

Denote by $D^b(\C_{X^{\rm an}})$ the derived category of $\C_{X^{\rm an}}$-modules. The full subcategory of $D^b(\C_{X^{\rm an}})$ consisting of objects whose cohomology are algebraically constructible (resp. algebraically constructible with support on $Z$) is denote by $D_{c}^b(X)$ (resp. $D_c^{b,\,Z}(X)$). Let ${\rm Perv}(X)$ be the full subcategory of $D^b_c(X)$ consisting of  perverse sheaves (\cite[Definition 8.1.20]{HTT}).

The de Rham functor for $M\in D^b(D_X)$ is defined by
$$\DR_X(M):=\Omega_{X^{\rm an}}\otimes_{D_{X^{\rm an}}}^LM^{\rm an}.$$
By \cite[Theorem 7.2.2]{HTT}, the de Rham functor $\DR_X$ gives two equivalences of categories:
\begin{eqnarray*}
\DR_X\colon D^b_{rh}(D_X)\simeq D^b_c(X)\hbox{ and }\DR_X\colon {\rm Mod}_{rh}(D_X)\simeq{\rm Perv}(X),
\end{eqnarray*}
which are now called the Riemann-Hilbert correspondence.

If $f\colon X\to Y$ is a morphism of schemes, we write $(f^{\rm an})^{-1},(f^{\rm an})^!, Rf^{\rm an}_*, Rf^{\rm an}_!$ as $f^{-1},f^!,f_*,f_!$, respectively. For any schemes $X$ and $Y$, we have two bi-functors
\begin{align*}&\boxtimes\colon D_c^b(X)\times D_c^b(Y)\to D_c^b(X\times_\C Y),&&\sF\boxtimes\sG:=p_X^{-1}\sF\otimes_{\C_{X^{\rm an}\times Y^{\rm an} }}p_Y^{-1}\sG;\\
&\boxtimes\colon {\rm Perv}(X)\times {\rm Perv}(Y)\to {\rm Perv}(X\times_\C Y),&&\sF\boxtimes\sG:=p_X^{-1}\sF\otimes_{\C_{X^{\rm an}\times Y^{\rm an} }}p_Y^{-1}\sG,
\end{align*}
where $p_X\colon X\times_\C Y\to X$ and $p_Y\colon X\times_\C Y\to Y$ are the projections.

For any immersion $f\colon X\to X'$ of schemes, we also have the minimal extension functor
$$f_{!*}\colon {\rm Perv}(X)\to {\rm Perv}(Y).$$
Actually, for any perverse sheaf $\sF$ on $X$, $f_{!*}\sF={\rm im}({^pf}_!\sF\to{^p}f_*\sF)$ (\cite[Definition 8.2.2]{HTT}).

If $g:Y\to Y'$ is another immersion, then for any $\sF\in{\rm Perv}(X)$ and $\sG\in{\rm Perv}(Y)$, we have
$$(f\times g)_{!*}(\sF\boxtimes\sG)=f_{!*}\sF\boxtimes g_{!*}\sG\in{\rm Perv}(X'\times_\C Y').$$

\subsection{Filtration on $N_A(\beta)$}

To define a filtration on $N_A(\beta)$ by combinatorial properties of the convex polyhedral cone generated by the column vectors of $A$, we use the notation $\mathbb S A=\sum\limits\limits_{j=1}^N\mathbb Sa_j\subset\C^n$ for any subset $\mathbb S$ of $\C$. In this paper, $\mathbb S$ could be $\N$, $\Z$, $\mathbb Q$, $\mathbb R$, $\mathbb R_{\geq0}$ or $\C$.

\begin{defn}
Identify the matrix $A$ with the set of column vectors of $A$. A \emph{face} of $A$ is a subset $F$ of $A$ such that there exists $h\in{\rm Hom}_\Z(\Z A,\,\Z)$ with the properties that $h(A)\subset\N$ and $F=\ker(h)\cap A$. We write $F\prec A$ if $F$ is a face of $A$. A \emph{facet} of $A$ is a face of codimension 1.
\end{defn}

For any $F\prec A$, $x_F$ is a subset of $x_A$. We also define $\C[x_F]$, $D_F$, $\A^F$, $T_F$ and $ i_F\colon X_F={\rm Spec}\,S_F\to\A^F$. Let $d_F$ be the rank of the matrix $F$.
Fix notations by the commutative diagram
\[\xymatrix{X_F\ar[dd]^{ i_F}\ar[rr]^{\bar i_{F,\,A}}&&X_A\ar[dd]^{ i_A}\\
&T_F\ar[lu]_{\bar j_F}\ar[ld]_{j_F}\ar[rd]^{j_{F,\,A}}\ar[ru]^{\bar j_{F,\,A}}&\\
\A^F\ar[rr]^{i_{F,\,A}}&&\A^A,}\]
where $i_{F,\,A}\colon \A^F\to\A^A$ is the closed immersion defined by the ideal of $\C[x_A]$ generated by those $x_j$ such that $a_j\notin F$, and where $\bar j_F\colon T_F\to X_F$ is the morphism induced by the homomorphism $\C[t^{a_j}\,|\,a_j\in F]\hookrightarrow\C[t^{\pm a_j}\,|\,a_j\in F]$.

For any $F\prec A$ with $\beta\in\C F$, the free $\sO_{T_F}$-module $\sO_{T_F}\cdot t^{-\beta}$ of rank one generated by the symbol $t^{-\beta}$ is equipped with a $D_{T_F}$-module structure via the product rule. Denote this $D_{T_F}$-module $\sO_{T_F}\cdot t^{-\beta}$ by $\mathcal{O}_{T_F}^\beta$.

The matrix $A$ is called \textit{normal} if $\mathbb R_{\geq0}A\cap\Z A=\N A$. This is equivalent to saying that the toric variety $X_A$ is normal.

\begin{defn}\label{che}
For any $0\leq i\leq d_A$, denote by $W_i(A, \beta)$ the $D_A$-submodule of $N_A(\beta)=D_A\Big/D_AI_A+\sum\limits\limits_{i=1}^nD_A\Big(\sum\limits\limits_{j=1}^Na_{ij}\p_{j}x_j+\beta_i\Big)$ generated by the image of $\prod\limits_{j=1}^Nx_j^{v_j}$ in $N_A(\beta)$ for those $v=(v_1,\ldots,v_N)^{\rm t}\in\N^N$ such that $Av\in\N A\backslash\mathop{\bigcup}\limits_{\substack{F\prec
A\\d_A-d_F>i}}\N F$. We get a filtration
$$0\subset W_0(A, \beta)\subset\cdots\subset W_{d_A}(A, \beta)=N_A(\beta)$$
on $N_A(\beta)$ which coincides with the Fourier transform of the filtration on $M_A(\beta)$ considered in \cite{A2}.\end{defn}

\begin{defn}\label{20181}
(1) For any facet $F$ of $A$, there is a unique linear form $\ell_F\colon \C A\to\C$, called the primitive integral support function of $F$, with the properties that $\ell_F(A)\subseteq\N$, $\ell_F(F)=0$ and $\ell_F(\Z A)=\Z$.

(2) For any $\beta\in\C A$, we say that $\beta$ is \textit{$A$-nonresonant (resp. weakly $A$-nonresonant,  resp. semi $A$-nonresonant)} if $\ell_F(\beta)\notin\Z$ (resp. $\ell_F(\beta)\notin\Z-\{0\}$, resp. $\ell_F(\beta)\notin\Z_{<0}$ ) for any facet $F$ of $A$.

(3) Let $F_1,\ldots,F_k$ be all facets of $A$ for which $\ell_{F_i}(\beta)\in \mathbb{Z}$. We say that $A$ is \emph{simplicial relative to $\beta$} if for all $i$, the intersection of any $i$ distinct elements of $\{F_1,\ldots,F_k\}$ is a face of codimension $i$ in $A$.
\end{defn}

In this paper, we prove the following theorem.

\begin{thm}\label{baaa}
(1) For any $0\leq i\leq d_A$, we have a canonical epimorphism
$$W_i(A, \beta)/W_{i-1}(A, \beta)\to\bigoplus_{\substack{F\in \mathfrak F_i(A)\\\beta\in\C F}}(i_{F,\,A})_+W_0(F, \beta)$$
of regular holonomic $D$-modules on $\A^A$, where $\mathfrak F_i(A)$ is the set of faces of $A$ of codimension $i$.

(2) If $A$ is simplicial relative to $\beta$, then the above epimorphism is an isomorphism.

(3) Let $F$ be a face of $A$ with $\beta\in\C F$. If $F$ is normal and $\beta$ is weakly $F$-nonresonant, then $W_0(F,\beta)=(j_F)_{!+}\mathcal{O}_{T_F}^\beta$. In particular, $W_0(F,\beta)$ is irreducible.

(4) Suppose that $A$ is simplicial relative to $\beta$ and that $\beta$ is weakly $A$-nonresonant. Given $0\leq i\leq d_A$ such that $F$ is normal for any $F\in\mathfrak F_i(A)$ with $\beta\in\C F$, then we have
$$W_i(A, \beta)/W_{i-1}(A, \beta)\simeq\bigoplus_{\substack{F\in \mathfrak F_i(A)\\\beta\in\C F}}(j_{F,\,A})_{!+}\mathcal{O}_{T_F}^\beta.$$
In particular, $W_i(A,\beta)/W_{i-1}(A,\beta)$ is semisimple.
\end{thm}

\begin{rem}
Part (3) of Theorem \ref{baaa} implies the irreducibility of $W_0(A,\,\beta)$ conjectured in \cite{A2} when $A$ is normal. In Example \ref{hou}, we show that $W_0(A,0)$ may not be semisimple if $A$ is not normal.
\end{rem}

\subsection{Filtration on perverse sheaves on $X_A$}

For any scheme $X$, we call a locally free $\C_{X^{\rm an}}$-module of finite rank a local system on $X$. Given a morphism $f:X\to Y$ of schemes and a local system $\sL$ on $X$, denote by $f^{\diamond}(\sL )$ the set of isomorphism classes of local systems on $Y$ whose inverse images on $X$ are isomorphic to $\sL$.

Let $\sL_A$ be a rank one local system on $T_A$. Then $\sL_A[d_A]$ is an irreducible perverse sheaf on $T_A$. The torus embedding $\bar j_A\colon T_A\to X_A$ is affine, $(\bar j_A)_*(\sL_A[d_A])$ is therefore a perverse sheaf on $X_A$. For any $0\leq i\leq d_A$, $\bar j_A$ factors as $T_A\x{\bar \ell_i^A}U_i(A)\x{\bar k_i^A}X_A,$
where
$U_i(A)=X_A-\mathop{\bigcup}\limits_{\substack{F\prec A\\d_F<d_A-i}}X_F.$
 Any face $F$ of $A$ defines a homomorphism $\pi_{F,A}:T_A\to T_F$ of tori.
\begin{defn}\label{18321}
(1) Define a filtration
$$0\subset\sW_0(\sL_A)\subset\cdots\subset\sW_{d_A}(\sL_A)=(\bar j_A)_*(\sL_A[d_A])$$
of the perverse sheaf $(\bar j_A)_*(\sL_A[d_A])$ on $X_A$ by perverse subsheaves
$\sW_i(\sL_A)=(\bar k_i^A)_{!*}(\bar \ell_i^A)_*(\sL_A[d_A]).$

(2) Let $F_1,\ldots,F_k$ be all facets of $A$ for which $\pi_{F_i,A}^\diamond(\sL_A)\neq\emptyset$. We say that $A$ is \emph{simplicial relative to $\sL_A$} if for all $i$, the intersection of any $i$ distinct elements of $\{F_1,\ldots,F_k\}$ is a face of codimension $i$ in $A$.
\end{defn}

The theorem corresponding to Theorem \ref{baaa} for perverse sheaves is the following.

\begin{thm}\label{yaode}
Let $\sL_A$ be a rank one local system on $T_A$.

(1) For any $0\leq i\leq d_A$, we have a canonical epimorphism
$$\sW_i(\sL_A)/\sW_{i-1}(\sL_A)\to\bigoplus_{\substack{F\in \mathfrak F_i(A)\\\sL_F\in\pi_{F,A}^{\diamond}(\sL_A)}}(\bar j_{F,\,A})_{!*}(\sL_F[d_F])$$
of perverse sheaves on $X_A$.

(2) If $A$ is simplicial relative to $\sL_A$, then the above epimorphism is an isomorphism.
\end{thm}

The relation between Theorem \ref{baaa} and Theorem \ref{yaode} is the following.

\begin{thm}\label{maaa}
Suppose that  $\beta\in\C A$ and let $0\leq i\leq d_A$. There exists a unique rank one local system $\sL_A$ on $T_A$ such that $\sL_A[d_A]={\rm DR}_{T_A}(\mathcal{O}_{T_A}^\beta)$.

(1) If $A$ is normal, then $\beta$ is semi $A$-nonresonant if and only if ${\rm DR}_{\A^A}(N_A(\beta))=(j_A)_*(\sL_A[d_A])$.

(2) If $A$ is normal and $\beta$ is weakly $A$-nonresonant, then ${\rm DR}_{\A^A}(W_i(A, \beta))\simeq( i_A)_*\sW_i(\sL_A)$.


(3) The parameter $\beta$ is $A$-nonresonant if and only if $N_A(\beta)$ is irreducible. In this case,
$$N_A(\beta)=(j_A)_{+}\sO_{T_A}^\beta=(j_A)_{!+}\sO_{T_A}^\beta\hbox{ and }{\rm DR}_{\A^A}(N_A(\beta))\simeq(j_A)_{*}(\sL_A[d_A])\simeq(j_A)_{!*}(\sL_A[d_A]).$$

\end{thm}

The paper is organized as follows. In section 2, we redefine GKZ-systems and their Fourier transforms by $D$-module functors. In section 3, we use the Euler-Koszul complex to study GKZ-systems.  In section 4, we study filtrations on $N_A(\beta)$ and $(j_A)_+\mathcal{O}_{T_A}^\beta$, respectively. In section 5, we study filtrations on the corresponding perverse sheaf $(\bar j_A)_*(\sL_A[d_A])$ on $X_A$.

\textbf{Acknowledgements.}
The author express his gratitude to Lei Fu for his suggestion and discussion on this project.
In addition, we would like to thank the referee for particularly careful readings, detailed and constructive comments which are very helpful to the improvement of this paper.
My research is supported by the NSFC grant No. 11671269.

\section{GKZ-systems and Fourier transform}
In $\S 2.1$ we state two functorial properties of $D$-modules which allow us to reduce Theorem \ref{baaa}, Theorem \ref{yaode} and Example \ref{hou} to the case when $A$ is normal or simplicial. In $\S 2.2$ we give variations of the GKZ-systems and their Fourier transforms by $D$-module functors.
\subsection{Exact $D$-module functors}
\begin{lem}\label{4131}
(1) Given a morphism $f\colon X\to Y$ of smooth schemes, let $Z$ be a closed subscheme of $X$, and let $T$ be a closed subscheme of $Y$ such that $f(Z)\subseteq T$. Suppose that $f|_Z\colon Z\to T$ is a finite morphism or an affine immersion.
Then $f_+M,\,f_!M\in{\rm Mod}_{rh}^{T}(D_Y)$ for any $M\in {\rm Mod}^Z_{rh}(D_X)$. We have two exact functors $f_+\hbox{ and }f_!\colon{\rm Mod}_{rh}^Z(D_X)\to{\rm Mod}_{rh}^T(D_Y)$, and they coincide if $f|_Z$ is finite.

(2) Consider a commutative diagram
\[\xymatrix{X_1\ar[r]^j\ar[d]^\tau&X_2\ar[d]^{\pi}\\X_3\ar[r]^k&X_4}\]
of smooth schemes, where $j$ and $k$ are immersions. Given a closed subscheme $Z_i$ of each $X_i$, assume that $j(Z_1)\subseteq Z_2,\,\tau(Z_1)\subseteq Z_3,\,\pi(Z_2)\subseteq Z_4$ and $k(Z_3)\subseteq Z_4$. If $\tau|_{Z_1}\colon Z_1\to Z_3$ and $\pi|_{Z_2}\colon Z_2\to Z_4$ are finite morphisms, then for any $M\in{\rm Mod}_{rh}^{Z_1}(D_{X_1})$, we have
$$\pi_+(j_{!+}M)=k_{!+}(\tau_+M)\in{\rm Mod}_{rh}^{Z_4}(D_{X_4}).$$ \end{lem}

\begin{proof}
(1) By the Riemann-Hilbert correspondence, we have a commutative diagram
\begin{eqnarray*}\xymatrix{D_{rh}^{b,\,Z}(D_X)\ar[d]^{f_+}\ar[r]^{{\rm DR}_X}&D_c^{b,\,Z}(X)\ar[d]^{f_*}\ar[r]^\simeq&D_c^b(Z)\ar[d]^{(f|_Z)_*}\\
D_{rh}^{b,\,T}(D_Y)\ar[r]^{{\rm DR}_Y}&D_c^{b,\,T}(Y)\ar[r]^\simeq&D_c^b(T)}
\end{eqnarray*}
of derived categories, whose horizontal arrows are equivalences of categories. The same holds if we replace $f_+$, $f_*$ and $(f|_Z)_*$ by $f_!$, $f_!$ and $(f|_Z)_!$, respectively. So part (1) follows from the corresponding properties for perverse sheaves.

For part (2), note that $j|_{Z_1}\colon Z_1\to Z_2$ and $k|_{Z_3}\colon Z_3\to Z_4$ are immersions. For simplicity, set $\bar j=j|_{Z_1}$, $\bar k=k|_{Z_3}$, $\bar \tau=\tau|_{Z_1}$ and $\bar\pi=\pi|_{Z_2}$. According to the proof of part (1), to prove part (2), we only need to verify that $\bar\pi_*(\bar j_{!*}\sF)=\bar k_{!*}(\bar \tau_*\sF)\in{\rm Perv}(Z_4)$ for any $\sF\in{\rm Perv}(Z_1)$. Since $\bar\pi_*=\bar\pi_!:D_c^b(Z_2)\to D_c^b(Z_4)$ and $\bar \tau_!=\bar \tau_*:D_c^b(Z_1)\to D_c^b(Z_3)$ are $t$-exact functors, we have $\bar\pi_*\bar j_!\sF=\bar\pi_!\bar j_!\sF=\bar k_!\bar \tau_!\sF=\bar k_!\bar \tau_*\sF\in D_c^{b}(Z_4)$ and hence
$$\bar\pi_*{^p}\bar j_!\sF={^pH^0}(\bar\pi_*\bar j_!\sF)={{^p}H}^0(\bar k_!\bar \tau_*\sF)={^p\bar k_!}\bar\tau_*\sF\in{\rm Perv}(Z_4).$$
Similarly, $\bar\pi_*{^p}\bar j_*\sF={^p\bar k_*}\bar\tau_*\sF\in{\rm Perv}(Z_4).$ It follows immediately that
$$\bar\pi_*\bar j_{!*}\sF=\bar\pi_*{\rm im}({^p}\bar j_!\sF\to{^p}\bar j_*\sF)={\rm im}(\bar\pi_*{^p}\bar j_!\sF\to\bar\pi_*{^p}\bar j_*\sF)={\rm im}({^p\bar k_!}\bar\tau_*\sF\to{^p\bar k_*}\bar\tau_*\sF)=\bar k_{!*}\bar \tau_*\sF.$$
\end{proof}

\begin{lem}\label{4132}
Let $B$ be a subset of $A$ such that $\mathbb R_{\geq0}A=\mathbb R_{\geq0}B$ or $A\subset B\cup(-\N B)$. The inclusion $B\subset A$ defines a morphism $\pi\colon \A^A\to\A^B$.
Then $\pi_+M\in{\rm Mod}^{X_B}_{rh}(D_{B})$ for any $M\in {\rm Mod}^{X_A}_{rh}(D_{A})$, and we have an exact functor $\pi_+\colon{\rm Mod}_{rh}^{X_A}(D_A)\to{\rm Mod}_{rh}^{X_B}(D_B)$. Moreover, this exact functor $\pi_+$ is faithful if $\mathbb R_{\geq0}A=\mathbb R_{\geq0}B$.
\end{lem}

\begin{proof}
If $\mathbb R_{\geq0}A=\mathbb R_{\geq0}B$, then by Exercise (iv) in \cite[Page 258]{BH}, the homomorphism $S_B\to S_A$ of finitely generated $\C$-algebras is integral. So $\pi|_{X_A}\colon X_A={\rm Spec}\,S_A\to X_B={\rm Spec}\,S_B$ is a finite surjective morphism. In this case, the lemma follows immediately from part (1) of Lemma \ref{4131}.

If $A\subset  B\cup(-\N B)$, then $A=B\cup\{-b_1,\ldots,-b_\ell\}$ for some $b_i\in\N B$.
Hence $S_A=S_B[t^{-b_1},\ldots,t^{-b_\ell}]$, so that $\pi|_{X_A}\colon X_A\to X_B$ is an open immersion of affine schemes. In this case, the lemma also follows from part (1) of Lemma \ref{4131}.
\end{proof}

\subsection{Variations of the GKZ-systems and their Fourier transforms} Let $\A^1={\rm Spec}\,\C[z]$. The free $\sO_{\A^1}$-module $\sO_{\A^1}\cdot e^{z}$ (resp. $\sO_{\A^1}\cdot e^{-z}$) of rank one generated by the symbol $e^{z}$ (resp. $e^{-z}$), is equipped with a $D_{\A^1}$-module structure via the product rule. We still denote this $D_{\A^1}$-module by $e^z$ (resp. $e^{-z}$).

\begin{defn}\label{xi}
For any $\beta\in\C A$, define
\begin{eqnarray}\label{1026a}
{\rm Hyp}_A(\beta)=p_{2+}(p_1^+ \mathcal{O}_{T_A}^\beta\otimes^L_{\sO_{T_A\times\A^A}} f_A^+ (e^z)),
\end{eqnarray}
where
$$p_1\colon T_A\times\A^A\to T_A,\;\;p_2\colon T_A\times\A^A\to\A^A$$
are the projection maps and
$$f_A\colon T_A\times\A^A\to\A^1$$ is the morphism defined by
$\sum\limits\limits_{j=1}^Nt^{a_j}\otimes x_j.$
\end{defn}


\begin{defn}
For any $M\in D^b(D_{A})$, define the Fourier transform of $M$ to be
$$\sF_A(M)=\pi_{2+}(\pi_1^+ M\otimes^L_{\sO_{\A^A\times\A^A}}\langle\,,\,\rangle^+ e^{-z})$$
where $\pi_i\colon \A^A\times\A^A\to\A^A$ is the projection onto the $i$-th factor and
$$\langle\,,\,\rangle\colon \A^A\times\A^A\to\A^1$$
is the morphism defined by $\sum\limits\limits_{j=1}^Nx_j\otimes x_j.$
\end{defn}

The Fourier transform $\sF_A$ preserves holonomicity and irreducibility but not regularity. We refer to \cite{Br} and \cite{DE} for more details of the Fourier transform of $D$-modules.

Suppose that  $\beta\in\C A$. Using the same method of \cite [Lemma 1.1] {F}, we have
 $$\sF_A({\rm Hyp}_A(\beta))\simeq (j_A)_+\mathcal{O}_{T_A}^\beta.$$
By \cite[4.5]{GKZ}, there is a canonical homomorphism $N_A(\beta)\to (j_A)_+\mathcal{O}_{T_A}^\beta$, and I will determine when $N_A(\beta)\simeq(j_A)_+\mathcal{O}_{T_A}^\beta$ in Lemma \ref{fei1}.
So the same assertion holds for $M_A(\beta)$ and ${\rm Hyp}_A(\beta)$ by the exactness of the Fourier transform.

\section{Euler-Koszul homology of toric modules}
In $\S 3.1$, we recall the definitions of Euler-Koszul complexes and toric modules given in \cite{MMW}, and in $\S 3.2$ we study their functorial and vanishing properties. In $\S 3.3$, we formulate some classes of parameters which are very important for filtrations on $N_A(\beta)$.

\subsection{Euler-Koszul complexes}
 L.F. Matusevich, E. Miller and U. Walther in \cite{MMW} introduced the Euler-Koszul complex to study the GKZ-system $M_A(\beta)$. For the purpose of this paper, we need to consider the Fourier transform of the Euler-Koszul complex in \cite{MMW}, which we still called the Euler-Koszul complex for simplicity. The Fourier transform $N_A(\beta)$ of the GKZ-system  $M_A(\beta)$ is the 0-th homology of the Euler-Koszul complex of the graded algebra $S_A$ (see Lemma \ref{1029c}), and the filtration of $S_A$ by homogenous ideals thus gives a filtration of $N_A(\beta)$ (see Definition \ref{2018}).

In this section, set $\epsilon_A=\sum\limits\limits_{j=1}^Na_j$. For any $F\prec A$, the rings $\C[x_F]$, $S_F$ and $D_F$ are naturally graded by $\Z^n$ by setting $-\deg(\p_{j})=\deg(x_j)=\deg(t^{a_j})=a_j$. For any $1\leq i\leq n$, let
$$E^A_i=\sum\limits_{j=1}^Na_{ij}\p_{j}x_j\in D_A.$$
The differential operator $E_i^A+\beta_i$ is homogenous of degree 0. For any homogenous differential operator $P\in D_A$ of degree $(d_1,\ldots,d_n)^{\rm t}$, we have
$[E_i^A+\beta_i,\,P]=d_iP.$ In particular, the $n$ operators $E_1^A+\beta_1,\ldots,E_n^A+\beta_n$ commute with each other.

\begin{defn}\label{qian}
Let $M$ be a $\Z^n$-graded $\C[x_A]$-module. Then $D_A\otimes_{\C[x_A]}M=\C[\p_{A}]\otimes_\C M$ is a $\Z^n$-graded $D_A$-module.

(1) Define the action of $E^A_i+\beta_i$ on $\C[\p_{A}]\otimes_\C M$ by
$$(E^A_i+\beta_i)(m)=(E^A_i+\beta_i-d_i).m$$
for any homogenous element $m\in \C[\p_{A}]\otimes_\C M $ of degree $(d_1,\ldots,d_n)^{\rm t}\in\Z^n.$ Then the action of $E_i^A+\beta_i$ on $\C[\p_A]\otimes_\C M$ is $D_A$-linear.

(2) \cite[Definition 4.2]{MMW} Write $E^A+\beta=\{E_1^A+\beta_1,\ldots,E_n^A+\beta_n\}$. The \emph{Euler-Koszul complex} $\sK_\bullet(E^A+\beta,\; M)$ is the Koszul complex of $D_A$-modules defined by the $n$ commutative operators $E^A_1+\beta_1,\ldots,E^A_n+\beta_n$ on the $D_A$-module $\C[\p_{A}]\otimes_\C M$. That is, $\sK_\bullet(E^A+\beta,\;M)=\sK(E^A+\beta,\;\C[\p_A]\otimes_\C M)$. The $i$-th homology of $\sK_\bullet(E^A+\beta,\;M)$ is denoted by $H_i(E^A+\beta,\;M)$. Then $H_i(E^A+\beta,\;M)=0$ unless $0\leq i\leq n$.
\end{defn}

\noindent{\bf Example:} Let $B=(1)$ and let $A$ be the empty face of $B$. Then $S_B=\C[x]$ is $\Z$-graded by $\deg(x^i)=i$ and the action of the differential operator $E^A=0\in D_A=\C$ on $\C[\p_A]\otimes_\C S_B=\C[x]$ is given by
$$E^A(x^i)=E^A. x^i-ix^i=-ix^i.$$
So $\sK_\bullet(E^A,S_B)$ is quasi-isomorphic to the complex $\Big(\C[x]\x{}\C[x],\;\sum\limits_i c_ix^i\mapsto-\sum\limits_i ic_ix^i \Big)$. This proves that $H_i(E^A,S_B)=\C$ for $i=0,1$.

\begin{defn}\label{1222}(1) \cite[Definition 4.5]{MMW} A $\Z^n$-graded $\C[x_A]$-module $M$ is called toric if it has a finite filtration
$$0=M_0\subset M_1\subset\cdots\subset M_\ell=M$$
by $\Z^n$-graded $\C[x_A]$-submodules such that for any $i$, $M_i/M_{i-1}$ is isomorphic to $S_F[b]$ as $\Z^n$-graded $\C[x_A]$ modules for some $F\prec A$ and $b\in\Z^n$.

(2) \cite[Definition 5.2]{MMW} For any toric $\C[x_A]$-module $M=\bigoplus\limits_{a\in\Z^n}M_a$, let $$\deg(M)=\{a\in\Z^n|\;M_a\neq 0\},$$ and let ${\rm qdeg}(M)$ be the Zariski closure of $\deg(M)$ in $\C^n$. Then ${\rm qdeg}(M)$ is a finite union of subsets of $\C^n$ of the form $a+\C  F$ for some $a\in\Z^n$ and $F\prec A$ such that $a+\N F\subset\deg(M)$.
\end{defn}

\subsection{Functorial and vanishing properties of the Euler-Koszul complexes}
Lemma \ref{1029c} and Lemma \ref{xiu} prove some functorial properties of Euler-Koszul complexes for closed immersions and projections, respectively.  Lemma \ref{lao} shows the regular holonomicity of $\sK_\bullet(E^A+\beta,S_A)$, which ensures that we can use perverse sheaves to study the GKZ-systems.  Corollary \ref{ahui} will be needed in the proof of Theorem \ref{baaa} and Example \ref{hou}.
\begin{lem}\label{1029c}
For any subset $F$ of $A$ and any $\Z^n$-graded $\C[x_F]$-module $M$, $M$ is also a $\Z^n$-graded $\C[x_A]$-module via the epimorphism $\C[x_A]\to\C[x_F]$
sending $x_j$ to $x_j$ if $a_j\in F$, and to $0$ otherwise.
This epimorphism defines a closed immersion $i_{F,\,A}\colon \A^F\to \A^A$. For any $i\in\Z$, we have
\begin{eqnarray*}
&&N_F(\beta)=H_0(E^F+\beta,\,S_F);\\&&(i_{F,\,A})_+\sK_\bullet(E^F+\beta,\,M)=\sK_\bullet(E^A+\beta,\,M);\\&&(i_{F,\,A})_+H_i(E^F+\beta,\,M)=H_i(E^A+\beta,\,M).
\end{eqnarray*}
\end{lem}

\begin{proof}
The first equality follows immediately from the definition of $N_F(\beta)$.
The last two equalities hold by the exactness of the functor $(i_{F,A})_+$ and
\begin{eqnarray*}
&&(i_{F,A})_+\sK_\bullet(E^F+\beta,M)\\&=&(i_{F,A})_+\sK(E^F+\beta,\C[\p_F]\otimes_\C M)\\
&=&\C[\p_{A\backslash F}]\otimes_\C\sK(E^F+\beta,\C[\p_F]\otimes_\C M)\\
&=&\sK\Big({{\rm id}_{\C[\p_{A\backslash F}]}}\otimes (E^F+\beta),\;\C[\p_{A\backslash F}]\otimes_\C\C[\p_F]\otimes_\C M\Big)\\
&=&\sK(E^A+\beta,\;\C[\p_A]\otimes_\C M)\\
&=&\sK_\bullet(E^A+\beta,\;M).
\end{eqnarray*}
\end{proof}

When $d_A=n$ and $\N A\cap-\N A=0$, Matusevich, Miller and Walther in \cite{MMW} proved that for any toric $\C[x_A]$-module $M$, $\beta\notin{\rm qdeg}(M)\iff H_i(E^A+\beta, M)=0\hbox{ for any }i\iff H_0(E^A+\beta,M)=0$. They also showed that each homology $H_i(E^A+\beta,M)$ is a holonomic $D_A$-module. In the following, we prove these results with no hypotheses on $A$.

\begin{lem}\label{lao}
Let $M$ be a toric $\C[x_A]$-module and $\beta\in\C^n$.

(1) Then $\beta\notin{\rm qdeg}(M)$ if and only if $H_0(E^A+\beta,M)=0$, if and only if $H_i(E^A+\beta,M)=0$ for any $i\in\Z$.


(2) We have $\sK_\bullet(E^A+\beta,\,M)\in D_{rh}^{b,\,X_A}(D_A)$.
\end{lem}

\begin{proof}
First note that any short exact sequence $0\to M'\to M\to M''\to0$ of toric $\C[x_A]$-modules gives a long exact sequence of $D_A$-modules:
$$\cdots\to H_i(E^A+\beta,\,M')\to H_i(E^A+\beta,\,M)\to H_i(E^A+\beta,\,M'')\to H_{i-1}(E^A+\beta,M')\to\cdots.$$
This implies that, if the lemma holds for $M'$ and $M''$, it also holds for $M$. So we can assume that $M=S_F[b]$ for some $F\prec A$ and $b\in\Z^n$ by the definition of toric modules. By Lemma \ref{1029c},
\begin{eqnarray*}\label{li}
\sK_\bullet(E^A+\beta,\,S_F[b])=\sK_\bullet(E^A+\beta+b,\,S_F)=(i_{F,\,A})_+\sK_\bullet(E^F+\beta+b,\,S_F).
\end{eqnarray*}
We can therefore assume that $M=S_A$.

If $\beta\notin\C A={\rm qdeg}(S_A)$, then the $\C$-linear span of $E^A+\beta$ contains a nonzero scalar, and hence $H_i(E^A+\beta,S_A)=0$ for each $i$. So, to prove the lemma, it suffices to show that $H_0(E^A+\beta,S_A)\neq0$ and $\sK_\bullet(E^A+\beta,S_A)\in D_{rh}^{X_A}(D_A)$ if $\beta\in\C A$. We prove them by induction on $d_A$ as follows.

From now on, assume that $\beta\in\C A$. Let $A'$ be a submatrix of $A$ composed of $d_A$ many linearly independent rows of $A$, and let $\beta'$ be the part of $\beta$ in the corresponding rows. Without loss of generality, assume that $A'$ is composed of the first $d_A$ rows of $A$. Then $S_A\simeq S_{A'}$ as $\C[x_A]$-modules, $\beta'\in\C A'$ and the $\C$-span of $E^A+\beta$ coincides with that of $E^{A'}+\beta'$. Thus for any $d_A<i\leq n$, the induced action of $E_i^A+\beta_i$ on each homology $H_j(E^{A'}+\beta',S_A)$ is the zero action. Thus we have
\begin{align*}
&\quad\;\sK_\bullet(E^A+\beta,S_A)\\
&=\sK(E^A_1+\beta_1,\ldots,E^A_n+\beta_n,\;\C[\p_A]\otimes_\C S_A)\\
&=\sK\big(E^A_{d_A+1}+\beta_{d_A+1},\ldots,E_n^A+\beta_n,\;\sK(E^A_1+\beta_1,\ldots,E^A_{d_A}+\beta_{d_A},\;\C[\p_A]\otimes_\C S_A)\big)\\
&=\sK(0,\ldots,0,\;\sK(E^A_1+\beta_1,\ldots,E^A_{d_A}+\beta_{d_A}\;,\C[\p_A]\otimes_\C S_A))\\
&=\sK(0,\ldots,0,\;\sK(E^{A'}_1+\beta_1',\ldots,E^{A'}_{d_{A'}}+\beta'_{d_{A'}},\;\C[\p_{A'}]\otimes_\C S_{A'})).
\end{align*}
Hence
\begin{align}\label{tian1}H_i(E^A+\beta,\,S_A)=\bigoplus_{i+d_A-n\leq j\leq i}H_j(E^{A'}+\beta',\;S_{A'})^{\binom{n-d_A}{i-j}}.\end{align}
In particular, $H_0(E^A+\beta,S_A)=H_0(E^{A'}+\beta',S_{A'})$, which is nonzero by \cite{A}. This proves part (1) for any $A$.

If $d_A=0$, then (\ref{tian1}) implies that $H_i(E^A+\beta,\;S_A)=(D_A/D_A\C[x_A])^{\binom ni}$, which is a regular holonomic $D_A$-module supported on $X_A={\rm Spec}\,\C$. This proves part (2) when $d_A=0$.

Now suppose $d_A>0$. According to (\ref{tian1}), $\sK_\bullet(E^A+\beta,S_A)\in D_{rh}^{b,\,X_A}(D_A)$ if and only if $\sK_\bullet(E^{A'}+\beta',S_{A'})\in D_{rh}^{b,\,X_{A'}}(D_{A'})=D_{rh}^{b,\,X_A}(D_A)$. So we can assume that $d_A=n$.

Recall from \cite[Corollary 3.9]{SW} that
$$\sK_\bullet(E^A+\beta,\;t^{-k\epsilon_A}S_A)=(j_A)_+\mathcal{O}_{T_A}^\beta\in{\rm Mod}_{rh}^{X_A}(D_A)\hbox{ for }k\gg0.$$
In order to prove $\sK_\bullet(E^A+\beta,S_A)\in D_{rh}^{b,\,X_A}(D_A)$, it suffices to show that $\sK_\bullet(E^A+\beta,\,t^{-k\epsilon_A}S_A/S_A)\in D_{rh}^{b,\,X_A}(D_A)$, according to
the distinguished triangle
$$\sK_\bullet(E^A+\beta,\;S_A)\to\sK_\bullet(E^A+\beta,\,t^{-k\epsilon_A}S_A)\to\sK_\bullet(E^A+\beta,\,t^{-k\epsilon_A}S_A/S_A) $$
in $D^b(D_A)$. Since
$$\dim{\rm qdeg}(t^{-k\epsilon_A}S_A/S_A)=\dim\Big(\bigcup_{i=1}^k-i\epsilon_A+{\rm qdeg}(S_A/t^{\epsilon_A}S_A)\Big)=\dim{\rm qdeg}(S_A/t^{\epsilon_A}S_A)<d_A,$$
then $t^{-k\epsilon_A}S_A/S_A$ has a filtration with successive quotients of the form $S_F[b]$ for some proper faces $F$ of $A$ and $b\in\Z^n$. By induction hypothesis to these $F$, $\sK_\bullet(E^F+\beta,\,S_F[b])=\sK_\bullet(E^F+\beta+b,\,S_F)\in D_{rh}^{b,\,X_F}(D_F)$, and hence $\sK_\bullet(E^A+\beta,\,S_F[b])\in D_{rh}^{b,\,X_A}(D_A)$ by Lemma \ref{1029c}. It follows immediately that $\sK_\bullet(E^A+\beta,S_A)\in  D_{rh}^{b,\,X_A}(D_A)$. This completes the proof of part (2).
\end{proof}

\begin{lem}\label{xiu}
Let $B$ be a subset of $ A$, and let $M$ be a $\Z^n$-graded $\C[x_{A}]$-module. The inclusion $B\subset A$ defines a morphism $\pi\colon \A^{A}\to\A^B$ and a $\Z^n$-graded $\C[x_B]$-module structure on $M$. Then
$$\pi_+\sK_\bullet(E^{ A}+\beta,\;M)=\sK_\bullet(E^B+\beta,\;M).$$
\end{lem}
\begin{proof}
We refer to \cite[Proposition 1.5.28]{HTT} for the definition of $\pi_+$ of the projection $\pi\colon \A^A=\A^{A\backslash B}\times\A^B\to\A^B$. On $\C[\p_A]\otimes_\C M$, there are two sets of operators $E^A+\beta=\{E^A_1+\beta_1,\ldots,E^A_n+\beta_n\}$ and $\p_{A\backslash B}=\{\p_j|\;a_j\in A\backslash B\}$, where the action of $\p_j$ is given by left mulitiplication. These operators commute with each other, so we have
\begin{eqnarray*}
\pi_+\sK_\bullet(E^{ A}+\beta,\;M)
&=&\pi_+\sK(E^A+\beta,\C[\p_A]\otimes_\C M)\\
&=&\sK(\p_{A\backslash B},\sK(E^A+\beta,\C[\p_A]\otimes_\C M))\\
&=&\sK(E^A+\beta,\sK(\p_{A\backslash B},\C[\p_A]\otimes_\C M))\\
&=&\sK(E^A+\beta,(\C[\p_A]/\p_{A\backslash B}\C[\p_A])\otimes_\C M)\\
&=&\sK(E^B+\beta,\C[\p_B]\otimes_\C M)\\
&=&\sK_\bullet(E^B+\beta,\;M),
\end{eqnarray*}
where the fifth equality holds by the following commutative diagram
\[\xymatrix{\C[\p_B]\otimes_\C M\ar[d]^\simeq\ar[r]^{E_i^B+\beta_i}&\C[\p_B]\otimes_\C
M\ar[d]^\simeq\\
(\C[\p_A]/\p_{A\backslash B}\C[\p_A])\otimes_\C M\ar[r]^{E_i^A+\beta_i}&(\C[\p_A]/\p_{A\backslash B}\C[\p_A])\otimes_\C M}\]
induced by the composition $\C[\p_B]\hookrightarrow\C[\p_A]\twoheadrightarrow\C[\p_A]/\p_{A\backslash B}\C[\p_A]$.
\end{proof}

\begin{cor}\label{ahui}
Let $B$ be a subset of $ A$ with $\mathbb R_{\geq0}A=\mathbb R_{\geq0}B$ or $A\subset B\cup(-\N B)$, and let $\pi:\A^A\to\A^B$ be the morphism induced by $B\subset A$. For any toric $\C[x_A]$-module $M$ and any $i$, we have $$\pi_+H_i(E^{ A}+\beta,\;M)\simeq H_i(E^B+\beta,\;M)\in{\rm Mod}_{rh}^{X_B}(D_B).$$
In particular,
\begin{align*}
\pi_+{\rm im}\Big(H_0(E^A+\beta,\,I_i(A))\to H_0(E^A+\beta,\,S_A)\Big)={\rm im}\Big(H_0(E^B+\beta,\,I_i(A))\to H_0(E^B+\beta,\,S_A)\Big).
\end{align*}

\end{cor}

\begin{proof}
Since $M$ is a toric $\C[x_A]$-module, then by Lemma \ref{lao},  $\sK_\bullet(E^{A}+\beta,\;M)\in{D}_{rh}^{b,\,X_{A}}(D_{A})$ and $H_i(E^{A}+\beta,\;M)\in{\rm Mod}_{rh}^{X_A}(D_A)$.  According to Lemma \ref{4132}, we have an exact functor $\pi_+\colon{\rm Mod}_{rh}^{X_{A}}(D_{A})\to {\rm Mod}_{rh}^{X_{B}}(D_{ B})$. It follows from Lemma \ref{xiu} and the exactness of $\pi_+$ that
\begin{align*}\begin{aligned}
&\quad\;\pi_+H_i(E^A+\beta,M)=\pi_+H_i(\sK_\bullet(E^A+\beta,M))\\&=H_i(\pi_+\sK_\bullet(E^A+\beta,M))=H_i(\sK_\bullet(E^B+\beta,\,M))=H_i(E^B+\beta,\,M)\in{\rm Mod}_{rh}^{X_B}(D_B).
\end{aligned}\end{align*}
In particular,  for the toric $\C[x_A]$-modules $I_i(A)$ and $S_A$, we have
\begin{align*}
\pi_+H_0(E^A+\beta,I_i(A))=H_0(E^B+\beta,\,I_i(A))\hbox{ and }
\pi_+H_0(E^A+\beta,S_A)=H_0(E^B+\beta,S_A).
\end{align*}
So the exact functor $\pi_+\colon{\rm Mod}_{rh}^{X_{A}}(D_{A})\to {\rm Mod}_{rh}^{X_{B}}(D_{ B})$  implies that
\begin{align*}\begin{aligned}
&\quad\;\pi_+{\rm im}\Big(H_0(E^{A}+\beta,\;I_i(A))\to H_0(E^{A}+\beta,\;S_{A})\Big)\\
&={\rm im}\Big(\pi_+H_0(E^{A}+\beta,\;I_i(A))\to \pi_+H_0(E^{A}+\beta,\;S_{A})\Big)\\
&={\rm im}\Big(H_0(E^B+\beta,\;I_i(A))\to H_0(E^B+\beta,\;S_{A})\Big).
\end{aligned}\end{align*}
\end{proof}

\subsection{Classes of parameters} In Definition \ref{2018}, we formulate some classes of parameters which are very important for filtrations on $N_A(\beta)$. Lemma \ref{xu} shows that Definition \ref{2018} coincides with Definition \ref{20181} when $A$ is normal, and in Example \ref{zy} we give two counterexamples to Lemma \ref{xu} without the normality condition.

\begin{defn}\label{2018}
For any integer $i$ and $F\prec A$, let $$I_i(F)=\bigoplus_{a\in\N F\backslash\mathop{\bigcup}\limits_{\substack{G\prec F\\d_F-d_G>i}}\N G}\C t^a.$$
Define
\begin{eqnarray*}
&&{\rm Res}(A)=\Z A+\bigcup_{F\in \mathfrak F_1(A)}\C F;\\
&&{\rm sRes}(A)=-\Z_{\geq1}\epsilon_A+{\rm qdeg}(S_A/t^{\epsilon_A}S_A);\\
&&{\rm dRes}(A)=\bigcup_{\substack{0\leq i<d_A\\k\geq2}}{\rm qdeg}(I_i(A)/I_i(A)^k);\\
&&{\rm wRes}(A)={\rm sRes}(A)\cup{\rm dRes}(A).
\end{eqnarray*}
\end{defn}
Clearly, ${\rm sRes}(A)\subset{\rm wRes}(A)\subset{\rm Res}(A)$. The sets ${\rm Res}(A)$ and ${\rm sRes}(A)$ were originally defined in \cite[2.9]{GKZ} and \cite[Definition 3.4]{SW}, respectively. By the following lemma, ${\rm dRes}(A)\subseteq-{\rm sRes}(A)$, with equality if $A$ is normal. Roughly speaking, ``dRes" is the dual of ``sRes" and ``wRes" is weaker than ``Res".


\begin{lem}\label{xu}
Suppose that  $\beta\in\C A$.

(1) Then $\beta$ is $A$-nonresonant if and only if $\beta\notin{\rm Res}(A)$.

(2) If $\beta\notin{\rm sRes}(A)$, then $\beta$ is semi $A$-nonresonant. The converse holds if $A$ is normal.

(3) If $\ell_F(\beta)\notin\Z_{>0}$ for any facet $F$ of $A$, then $\beta\notin{\rm dRes}(A)$. The converse holds if $A$ is normal.

(4) Suppose that  $A$ is normal. Then $\beta\notin{\rm wRes}(A)$ if and only if $\beta$ is weakly $A$-nonresonant.
\end{lem}
\begin{proof}
Part (1) is trivial, and part (4) follows from (2) and (3).

For part (2), suppose that  $\ell_F(\beta)\in\Z_{<0}$ for some facet $F$ of $A$. Since $\ell_F(\Z A)=\Z$, there exists $m\in\N$ and $b\in\N A$ such that $\ell_F(\beta)=\ell_F(b-m\epsilon_A)$. Choose such a pair $(b,\,m)$ with $m$ as small as possible. Obviously, $m\geq1$.
For any $c\in \N F$, $\ell_F((b+c-\epsilon_A)-(m-1)\epsilon_A)=\ell_F(\beta)\in\Z_{<0}$. The minimality of $m$ implies that $b+c-\epsilon_A\notin\N A$ and then $b+c\in\deg(S_A/t^{\epsilon_A}S_A)$.
So $b+\N F\subset\deg(S_A/t^{\epsilon_A}S_A)$ and $\beta\in -m\epsilon_A+b+\C F\subset{\rm sRes}(A).$

Conversely, suppose that $\beta\in{\rm sRes}(A)$ and that $A$ is normal. By the definition of ${\rm sRes}(A)$ and considering the possible factor modules in a toric filtration of $S_A/t^{\epsilon_A}S_A$, $\beta=-m\epsilon_A+b+\alpha$ for some $m\in\Z_{>0}$, $b\in\N A$ and $\alpha\in\C^n$, such that there exists $F'\prec A$ with the properties that $\alpha\in\C F'$ and $b+\N F'\subset\deg(S_A/t^{\epsilon_A}S_A)$. For any facet $F$ of $A$ containing $F'$, we have
$$\ell_F(\beta)=\ell_F(-m\epsilon_A+b+\alpha)=\ell_F(b-\epsilon_A)-(m-1)\ell_F(\epsilon_A)\in\Z.$$
Then the converse of part (2) follows immediately from the claim below.

\underline{Claim:} There exists a facet $F$ of $A$ containing $F'$ with $\ell_F(b -\epsilon_A)<0$.

We prove this claim by contradiction. Suppose $\ell_F(b-\epsilon_A)\geq0$ for all facets $F$ of $A$ containing $F'$. Then by normality, $b-\epsilon_A \in  \N A -\N F'$, so there exists $c\in\N F'$ with $b + c \in \epsilon_A + \N A$. Hence, $b + \N F'$ is not contained in $\deg(R_A/t^{\epsilon_A}R_A)$. This proves the claim.

For part (3), suppose that  $\beta\in{\rm qdeg}(I_i(A)/I_i(A)^k)$ for some $0\leq i<d_A$ and $k\geq2$. Then $\beta=b+\alpha$ for some $b\in\N A$ and $\alpha\in\C^n$ such that there exists $F_0\prec A$ with the properties that $\alpha\in\C F_0$ and $b+\N F_0\subset\deg(I_i(A)/I_i(A)^k)$. For any $c\in\N F_0$, ${b+kc}\notin \deg(I_i(A)^k)$ implies that $c\notin\deg(I_i(A))$. Therefore, $d_{F_0}<d_A-i$ and $b\notin \mathbb R_{\geq0} F_0$.  Consequently, there exists a facet $F$ of $A$ containing $F_0$ such that $b\notin \mathbb R_{\geq0}F$. This shows that $\ell_F(\beta)=\ell_F(b)+\ell_F(\alpha)=\ell_F(b)\in\Z_{>0}$.

Conversely, suppose that $A$ is normal and that $\ell_F(\beta)\in\Z_{>0}$ for some facet $F$ of $A$. By (2) of \cite[Lemma 5.1]{S}, $\ell_F(\N A)=\N$. Then there exists $a\in\N A\backslash \N F$ such that $\ell_F(\beta)=\ell_F(a)\in\Z_{>0}$.  For any $c\in\deg(I_0(A)^k)$, $\ell_F(c)\geq k$, which implies that $a+\epsilon_F+\N F\subset\deg(I_0(A)/I_0(A)^k)$ for any $k>\ell_F(a)$ where $\epsilon_F=\sum\limits_{a_j\in F}a_j$. This proves that $\beta\in a+\C F\subseteq{\rm qdeg}(I_0(A)/I_0(A)^k)\subseteq{\rm dRes}(A).$
\end{proof}

We give two counterexamples of the converse of parts (2) and (3) in Lemma \ref{xu} without the normality assumption on $A$.

\begin{exam}\label{zy}
According to Definition \ref{20181}, set
\begin{eqnarray*}
&&{\rm SRes}(A)=\{\beta\in\C A|\;\ell_F(\beta)\in\Z_{<0}\hbox{ for some facet } F\hbox{ of } A\};\\
&&{\rm DRes}(A)=\{\beta\in\C A|\;\ell_F(\beta)\in\Z_{>0}\hbox{ for some facet } F\hbox{ of } A\}.
\end{eqnarray*}

(1) For $A=(2,3)$, the sets ${\rm sRes}(A)$, ${\rm dRes}(A)$, ${\rm SRes}(A)$ and ${\rm DRes}(A)$ are sketched below.
\\
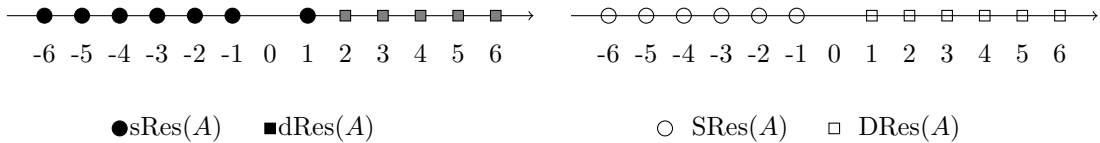
\begin{center}
\newdimen\scale
\scale=0.5cm
\begin{tikzpicture}
 \draw[->](-3.5,0)--(3.5,0);
 \foreach \x in {-6,...,-1,1}{
   \foreach \y in {0}{
     \node[draw,circle,inner sep=2pt,fill] at (\scale*\x,\scale*\y) {}; } }
\foreach \x in {2,...,6}{
   \foreach \y in {0}{
     \node[draw,rectangle,inner sep=2pt,fill=gray] at (\scale*\x,\scale*\y) {}; } }
\node[] at (0,-0.5) {0};\node[] at (0.5,-0.5) {1};\node[] at (1,-0.5) {2};\node[] at (1.5,-0.5) {3};\node[] at (2,-0.5) {4};\node[] at (2.5,-0.5) {5};\node[] at (3,-0.5) {6};\node[] at (-0.5,-0.5) {-1};\node[] at (-1,-0.5) {-2};\node[] at (-1.5,-0.5) {-3};\node[] at (-2,-0.5) {-4};\node[] at (-2.5,-0.5) {-5};\node[] at (-3,-0.5) {-6};
 \draw[->](4,0)--(11,0);
 \foreach \x in {9,...,14}{
   \foreach \y in {0}{
     \node[draw,circle,inner sep=2pt] at (\scale*\x,\scale*\y) {}; } }
\foreach \x in {16,...,21}{
   \foreach \y in {0}{
     \node[draw,rectangle,inner sep=2pt] at (\scale*\x,\scale*\y) {}; } }
\node[] at (7.5,-0.5) {0};\node[] at (7,-0.5) {-1};\node[] at (6.5,-0.5) {-2};\node[] at (6,-0.5) {-3};\node[] at (5.5,-0.5) {-4};\node[] at (5,-0.5) {-5};\node[] at (4.5,-0.5) {-6};\node[] at (8,-0.5) {1};\node[] at (8.5,-0.5) {2};\node[] at (9,-0.5) {3};\node[] at (9.5,-0.5) {4};\node[] at (10,-0.5) {5};\node[] at (10.5,-0.5) {6};
\node[draw,circle,inner sep=2pt,fill] at (-2,-1.5) {};
\node[] at (-1.25,-1.5) {${\rm sRes}(A)$};
\node[draw,rectangle,inner sep=2pt,fill] at (0,-1.5) {};
\node[] at (0.75,-1.5) {${\rm dRes}(A)$};
\node[draw,circle,inner sep=2pt] at (5.25,-1.5) {};
\node[] at (6.25,-1.5) {${\rm SRes}(A)$};
\node[draw,rectangle,inner sep=2pt] at (7.5,-1.5) {};
\node[] at (8.5,-1.5) {${\rm DRes}(A)$};
\end{tikzpicture}
\end{center}

(2) For
$A=\left(\begin{array}[c]{lll}1&1&0\\0&1&2\end{array}\right)$, the sets ${\rm sRes}(A)$, ${\rm dRes}(A)$, $\mathbb R_{\geq0}A$, $\N A$, ${\rm Res}(A)$, ${\rm SRes}(A)$ and ${\rm DRes}(A)$ are sketched below.

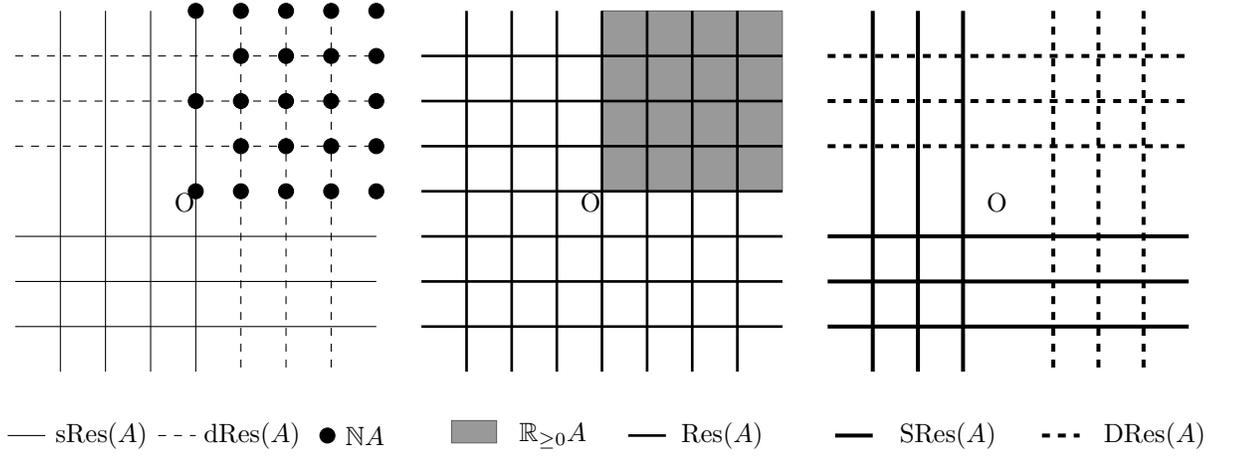
\begin{center}
\newdimen\scale
\scale=0.6cm
\begin{tikzpicture}
 \filldraw[black,opacity=.4] (\scale*9,\scale*0) -- (\scale*9,\scale*4) -- (\scale*13,\scale*4)--(\scale*13,\scale*0)  -- (\scale*9,\scale*0);
\draw[dashed] (\scale*-4,\scale* 3) -- (\scale*4, \scale *3);
\draw[dashed] (\scale*-4,\scale* 2) -- (\scale*4, \scale *2);
\draw[dashed] (\scale*-4,\scale* 1) -- (\scale*4, \scale *1);
\draw(\scale*-4,\scale* -1) -- (\scale*4, \scale *-1);
\draw (\scale*-4,\scale* -2) -- (\scale*4, \scale *-2);
\draw (\scale*-4,\scale* -3) -- (\scale*4, \scale *-3);
\draw[dashed] (\scale*3,\scale* 4) -- (\scale*3, \scale *-4);
\draw[dashed] (\scale*2,\scale* 4) -- (\scale*2, \scale *-4);
\draw[dashed] (\scale*1,\scale* 4) -- (\scale*1, \scale *-4);
\draw (\scale*0,\scale* 4) -- (\scale*0, \scale *-4);
\draw (\scale*-1,\scale* 4) -- (\scale*-1, \scale *-4);
\draw (\scale*-2,\scale* 4) -- (\scale*-2, \scale *-4);
\draw (\scale*-3,\scale* 4) -- (\scale*-3, \scale *-4);
\foreach \x in {1,...,4}{
   \foreach \y in {0,1,...,4}{
     \node[draw,circle,inner sep=2pt,fill] at (\scale*\x,\scale*\y) {}; } }
\foreach \x in {0}{
   \foreach \y in {0,2,4}{
     \node[draw,circle,inner sep=2pt,fill] at (\scale*\x,\scale*\y) {}; } }

\draw[line width=1pt](\scale*5,\scale* 3) -- (\scale*13, \scale *3);
\draw[line width=1pt] (\scale*5,\scale* 2) -- (\scale*13, \scale *2);
\draw[line width=1pt](\scale*5,\scale* 1) -- (\scale*13, \scale *1);
\draw[line width=1pt](\scale*5,\scale* 0) -- (\scale*13, \scale *0);
\draw[line width=1pt] (\scale*5,\scale* -1) -- (\scale*13, \scale *-1);
\draw[line width=1pt] (\scale*5,\scale* -2) -- (\scale*13, \scale *-2);
\draw[line width=1pt] (\scale*5,\scale* -3) -- (\scale*13, \scale *-3);
\draw[line width=1pt](\scale*12,\scale* 4) -- (\scale*12, \scale *-4);
\draw[line width=1pt] (\scale*11,\scale* 4) -- (\scale*11, \scale *-4);
\draw[line width=1pt](\scale*10,\scale* 4) -- (\scale*10, \scale *-4);
\draw[line width=1pt] (\scale*9,\scale* 4) -- (\scale*9, \scale *-4);
\draw[line width=1pt] (\scale*8,\scale* 4) -- (\scale*8, \scale *-4);
\draw[line width=1pt] (\scale*7,\scale* 4) -- (\scale*7, \scale *-4);
\draw[line width=1pt] (\scale*6,\scale* 4) -- (\scale*6, \scale *-4);

\draw[dashed, line width=1.5pt] (\scale*14,\scale* 3) -- (\scale*22, \scale *3);
\draw[dashed, line width=1.5pt] (\scale*14,\scale* 2) -- (\scale*22, \scale *2);
\draw[dashed, line width=1.5pt] (\scale*14,\scale* 1) -- (\scale*22, \scale *1);
\draw[line width=1.5pt] (\scale*14,\scale* -1) -- (\scale*22, \scale *-1);
\draw[line width=1.5pt] (\scale*14,\scale* -2) -- (\scale*22, \scale *-2);
\draw[line width=1.5pt] (\scale*14,\scale* -3) -- (\scale*22, \scale *-3);
\draw[dashed, line width=1.5pt]  (\scale*21,\scale* 4) -- (\scale*21, \scale *-4);
\draw[dashed, line width=1.5pt]  (\scale*20,\scale* 4) -- (\scale*20, \scale *-4);
\draw[dashed, line width=1.5pt] (\scale*19,\scale* 4) -- (\scale*19, \scale *-4);
\draw[line width=1.5pt] (\scale*17,\scale* 4) -- (\scale*17, \scale *-4);
\draw[line width=1.5pt] (\scale*16,\scale* 4) -- (\scale*16, \scale *-4);
\draw[line width=1.5pt] (\scale*15,\scale* 4) -- (\scale*15, \scale *-4);
 \foreach \x in {-0.25,8.75,17.75}{
   \foreach \y in {-0.25}{
     \node[] at (\scale*\x,\scale*\y) {O};}}
\draw(-2.5,-3.25)--(-2,-3.25);
\node[] at (-1.25,-3.25) {${\rm sRes}(A)$};
\draw[dashed](-0.5,-3.25)--(0,-3.25);
\node[] at (0.75,-3.25) {${\rm dRes}(A)$};
\node[draw,circle,inner sep=2pt,fill] at (1.75,-3.25) {};
\node[] at (2.25,-3.25) {$\N A$};
\filldraw[black,opacity=.4] (4,-3.05) -- (3.4,-3.05) -- (3.4,-3.35) -- (4,-3.35) -- (4,-3.05);
\node[] at (4.75,-3.25) {$\mathbb R_{\geq 0} A$};
\draw[line width=1pt](5.75,-3.25)--(6.25,-3.25);
\node[] at (7,-3.25) {${\rm Res}(A)$};
\draw[line width=1.5pt](8.5,-3.25)--(9,-3.25);
\node[] at (10,-3.25) {${\rm SRes}(A)$};
\draw[dashed, line width=1.5pt] (11.25,-3.25)--(11.75,-3.25);
\node[] at (12.75,-3.25) {${\rm DRes}(A)$};
\end{tikzpicture}
\end{center}
\end{exam}

\section{Toric filtrations on $D_A$-modules}
In this section, we mainly prove Theorem \ref{baaa}. In the process, we also prove Theorem \ref{maaa}, which follows immediately from Corollary \ref{fei} via the the Riemann--Hilbert correspondence. Corollary \ref{fei} also implies the irreducibility of $W_0(A,\,\beta)$ conjectured by Adolphson in \cite{A2} when $A$ is normal. If $A$ is not normal, we give a counterexample of this conjecture in Example \ref{hou}.

\subsection{Filtrations on $N_A(\beta)$ and $(j_A)_+\sO_{T_A}^\beta$} In this subsection, we define filtrations on $N_A(\beta)$ by Euler-Koszul complexes, and define filtrations on $(j_A)_+\sO_{T_A}^\beta$ by $D$-module functors.
The relation between these two filtrations will be discussed in Lemma \ref{fei1} and Corollary \ref{fei}. Lemma \ref{18317} will be used to calculate the composition factors of $N_A(0)$ in Example \ref{hou}.

Recall from Definition \ref{2018} that for any $F\prec A$, there is a filtration
$0\subset I_0(F)\subset\cdots\subset I_{d_F}(F)$
of $S_F=I_{d_F}(F)$ by ideals. Put
$$V_i(A)=\A^A-\bigcup_{\substack{F\prec A\\d_F<d_A-i}}\A^F,\;U_i(A)=X_A-\bigcup_{\substack{F\prec A\\d_F<d_A-i}}X_F\hbox{ and }Z_i(A)=X_A-U_i(A).$$
Thus $Z_i(A)={\rm Spec}\,S_A/I_i(A)$ and the morphism $j_A\colon T_A\to\A^A$ factors as $T_A\x{\ell_i^A}V_i(A)\x{k_i^A}\A^A.$
The immersions $j_A$ and $\ell_i^A$ are affine morphisms, then by part (1) of Lemma \ref{4131}, $(j_A)_+ \mathcal{O}_{T_A}^\beta\in{\rm Mod}_{rh}^{X_A}(D_A)$ and  $(\ell_i^A)_+ \mathcal{O}_{T_A}^\beta \in{\rm Mod}_{rh}^{U_i(A)}(D_{V_i(A)}).$

\begin{defn}\label{dimao}
(1) Define a filtration $$0\subset \widetilde W_0(A, \beta)\subset\cdots\subset\widetilde  W_{d_A}(A, \beta)=(j_A)_+\mathcal{O}_{T_A}^\beta $$ of $(j_A)_+ \mathcal{O}_{T_A}^\beta $ by $$\widetilde W_i(A, \beta)=(k_i^A)_{!+}(\ell_i^A)_+ \mathcal{O}_{T_A}^\beta .$$

(2) Define a filtration $$0\subset W_0(A, \beta)\subset\cdots\subset W_{d_A}(A, \beta)=N_A(\beta)$$ of $N_A(\beta)=H_0(E^A+\beta,\;S_A)$ by
$$W_i(A, \beta)={\rm im}\Big(H_0(E^A+\beta,\;I_i(A))\to H_0(E^A+\beta,\;S_A)\Big).$$
\end{defn}

The lemma below determines when the above two filtrations coincides.

\begin{lem}\label{fei1}
Suppose that  $\beta\in\C A$ and $0\leq i\leq d_A$.


(1) Then $\beta\notin{\rm sRes}(A)$ if and only if $\widetilde W_{d_A}(A,\beta)=W_{d_A}(A,\beta)$. In this case, $N_A(\beta)=(j_A)_+\sO_{T_A}^\beta$.

(2) If $\beta\notin{\rm wRes}(A)$, then $\widetilde W_i(A, \beta)=W_i(A, \beta)$. In this case, $W_0(A,\beta)=(j_A)_{!+}\sO_{T_A}^\beta$.

(3) Then $\beta\notin{\rm Res}(A)$ if and only if $N_A(\beta)$ is irreducible. In this case, $\widetilde W_i(A,\beta)=W_{i}(A,\beta)=(j_A)_{!+}\mathcal{O}_{T_A}^\beta$.
\end{lem}

\begin{proof}
By the same argument as in Lemma \ref{lao}, we can assume that $d_A=n$. Then (1) holds by \cite [Corollary 3.8]{SW}, and the only if part of (3) holds by \cite[Proposition 4.4, Theorem 4.6]{GKZ}.

For the if part of (3), suppose $\beta\in{\rm Res}(A)$ and prove by contradiction that $N_A(\beta)$ is reducible. Let $F$ be a minimal face of $A$ with $\beta\in\Z A+\C F$. By the definition of ${\rm Res}(A)$, $F$ is a proper face of $A$.
If $N_A(\beta)$ is irreducible, then so is $M_A(\beta)$. Consequently, $\C(x_A)\otimes_{\C[x_A]}M_A(\beta)$ is an irreducible $\C(x_A)[\p_A]$-module, where $\C(x_A)$ is the field of fractions of $\C[x_A]$. According to \cite[Lemma 3.5, Theorem 4.1]{SW1}, $\Z A=\Z F\oplus(\mathop{\bigoplus}\limits_{a_j\notin F}\Z a_j)$, and then $\beta=\beta'+\sum\limits_{a_j\notin F}n_ja_j$ for some $\beta'\in\C F$ and $n_j\in\Z$. So
$$N_A(\beta)=N_F(\beta')\otimes_\C\Big(\bigotimes_{a_j\notin F}\C[x_j,\p_j]/\C[x_j,\p_j](\p_jx_j+n_j)\Big).$$ Then the if part of (3) follows immediately from the reducibility of the left $\C[x_j,\p_j]$-module $\C[x_j,\p_j]/\C[x_j,\p_j](\p_j x_j+n_j)$.

To prove part (2), suppose that  $\beta\notin{\rm wRes}(A)$. By Definition \ref{2018}, $\beta\notin{\rm sRes}(A)$ and $\beta\notin{\rm qdeg}( I_i(A)/I_i(A)^k)$ for any $i$ and $k$. Part (1) implies that $N_A(\beta)=(j_A)_+\mathcal{O}_{T_A}^\beta$. Let $M_0$ be the $\C[x_A]$-submodule of $N_A(\beta)=D_A/D_AI_A+\sum\limits_{i=1}^nD_A(E^A_i+\beta_i)$ generated by the image of $1\in D_A$. So $(j_A)_+\mathcal{O}_{T_A}^\beta=N_A(\beta) $ is generated by $M_0$ as a $D_A$-module. Since $(j_A)_+ \mathcal{O}_{T_A}^\beta /\widetilde W_i(A, \beta)$ is supported on $(\A^A-V_i(A))\cap X_A=Z_i(A)$, then $J_i(A)^kM_0\subset\widetilde W_i(A, \beta)$ for $k\gg0$, where $J_i(A)$ is the ideal of $\C[x_A]$ defined by the closed subscheme $Z_i(A)$ of $\A^A$. By the definition of $(k_i^A)_{!+}$, $\widetilde W_i(A, \beta)=(k_i^A)_{!+}(\ell_i^A)_+\mathcal{O}_{T_A}^\beta $ is the smallest $D_A$-submodule of $(j_A)_+ \mathcal{O}_{T_A}^\beta $ whose restriction on $V_i(A)$ coincides with $(\ell_i^A)_+\mathcal{O}_{T_A}^\beta$. So
$$\widetilde W_i(A, \beta)=D_A.J_i(A)^kM_0={\rm im}\Big(H_0(E^A+\beta,\;J_i(A)^k )\to H_0(E^A+\beta,\; S_A)\Big)\hbox{ for }k\gg0.$$
Because the composition $J_i(A)^k\hookrightarrow\C[x_A]\twoheadrightarrow S_A$ also factors as $J_i(A)^k\twoheadrightarrow I_i(A)^k\hookrightarrow S_A$ and the functor $H_0(E^A+\beta,\,\bullet)$ is right exact, we get
\begin{eqnarray*}\label{1221}
\widetilde W_i(A, \beta)={\rm im}\Big(H_0(E^A+\beta,\;I_i(A)^k )\to H_0(E^A+\beta,\; S_A)\Big),\hbox{ for }k\gg0.
\end{eqnarray*}
According to $\beta\notin{\rm qdeg}(I_i(A)/I_i(A)^k)$ and (1) of Lemma \ref{lao}, we have $H_0(E^A+\beta,\;I_i(A)^k )=H_0(E^A+\beta,\;I_i(A))$, which implies that $\widetilde W_i(A,\beta)=W_i(A,\beta)$.



\end{proof}

Lemma \ref{xu} and Lemma \ref{fei1} have the following immediate corollary.

\begin{cor}\label{fei}
Let $\beta\in\C A$.

(1) If $A$ is normal, then $\beta$ is semi $A$-nonresonant if and only if $N_A(\beta)=(j_A)_+\mathcal{O}_{T_A}^\beta$.

(2) If $A$ is normal and $\beta$ is weakly $A$-nonresonant, then $\widetilde W_i(A, \beta)=W_i(A, \beta)$ for any $0\leq i\leq d_A$. In particular, $W_0(A,\beta)=(j_A)_{!+}\mathcal{O}_{T_A}^\beta$ is an irreducible regular holonomic $D$-module on $\A^A$.

(3) The parameter $\beta$ is $A$-nonresonant if and only if $N_A(\beta)$ is irreducible.
\end{cor}

In the following, we reformulate the filtration $\widetilde W_i(A,\beta)$ by the filtration $W_i(\widetilde A,\beta)$ for a normalization $\widetilde A$ of $A$.
\begin{lem}\label{18317}
Let $A$ be a submatrix of $\widetilde A$ of the same number of rows, and let $\pi\colon\A^{\widetilde A}\to\A^A$ be the induced morphism.
Suppose that $\beta\in\C A$  is weakly $A$-nonresonant and that $\mathbb R_{\geq0} A\cap\Z A=\N\widetilde A$. Then for any $0\leq i\leq d_A$, there is a canonical homomorphism $W_i(A, \beta)\to\widetilde W_i(A, \beta)$, and
\begin{eqnarray}\label{sun}\label{sun1}
&&\widetilde W_i(A, \beta)=\pi_+W_i(\widetilde A,\beta);\\\label{sun2}
&&\widetilde W_i(A, \beta)/\widetilde W_{i-1}(A,\beta)=\pi_+\big(W_i(\widetilde A,\beta)/W_{i-1}(\widetilde A,\beta)\big).
\end{eqnarray}
\end{lem}

\begin{proof}
The inclusion $A\subset\widetilde A$ defines a commutative diagram
\[\xymatrix{T_{\widetilde A}\ar[r]^{\ell_i^{\widetilde A}}\ar @{=} [d] &V_i(\widetilde A)\ar[r]^{k_i^{\widetilde A}}\ar[d]^{\pi_i}&\A^{\widetilde A}\ar[d]^{\pi}\\
T_A\ar[r]^{\ell_i^A}&V_i(A)\ar[r]^{k_i^A}&\A^A,}\]
whose horizontal morphisms are immersions. Applying part (1) of Lemma \ref{4131} to the affine immersion $\ell_i^{\widetilde A}$, we have $(\ell_i^{\widetilde A})_+\sO_{T_A}^\beta\in{\rm Mod}_{rh}^{U_i(\widetilde A)}(D_{V_i(\widetilde A)})$. As in the proof of Lemma \ref{4132}, the restriction morphisms $\pi|_{X_{\widetilde A}}\colon X_{\widetilde A}\to X_A$ and $\pi_i|_{U_i(\widetilde A)}\colon U_i(\widetilde A)\to U_i(A)$ are finite.  Applying part (2) of Lemma \ref{4131} to the right square of the above diagram and $(\ell_i^{\widetilde A})_+\sO_{T_A}^\beta\in{\rm Mod}_{rh}^{U_i(\widetilde A)}(D_{V_i(\widetilde A)})$, we thus have
\begin{align*}
(k_i^{A})_{!+}(\pi_i)_+(\ell_i^{\widetilde A})_{+}\mathcal{O}_{T_A}^\beta=\pi_+(k_i^{\widetilde A})_{!+}(\ell_i^{\widetilde A})_{+}\mathcal{O}_{T_A}^\beta= \pi_+\widetilde W_i(\widetilde A,\beta).
\end{align*}
It follows immediately that
\begin{eqnarray}\label{12231}\begin{split}
\widetilde W_i(A,\beta)=(k_i^A)_{!+}(\ell_i^A)_+\mathcal{O}_{T_A}^\beta=(k_i^{A})_{!+}(\pi_i)_+(\ell_i^{\widetilde A})_{+}\mathcal{O}_{T_A}^\beta=\pi_+\widetilde W_i(\widetilde A,\beta).
\end{split}\end{eqnarray}
Since $\beta$ is weakly $A$-nonresonant, it is weakly $\widetilde A$-nonresonant. Using part (2) of Corollary \ref{fei} with the normal matrix $\widetilde A$, we therefore have
\begin{align*}\widetilde W_i(\widetilde A,\beta)=W_i(\widetilde A,\beta).
\end{align*}
Combining this with (\ref{12231}), we prove (\ref{sun1}), and then by Lemma \ref{4132} we have
\begin{align*}
\widetilde W_i(A, \beta)/\widetilde W_{i-1}(A,\beta)=\pi_+W_i(\widetilde A,\beta)/\pi_+W_{i-1}(\widetilde A,\beta)=\pi_+\big(W_i(\widetilde A,\beta)/W_{i-1}(\widetilde A,\beta)\big).
\end{align*}

Applying Corollary \ref{ahui} to the inclusion $A\subset\widetilde A$, we get
\begin{align*}\begin{aligned}
\pi_+W_i(\widetilde A,\,\beta)={\rm im}\Big(H_0(E^A+\beta,\;I_i(\widetilde A))\to H_0(E^A+\beta,\;S_{\widetilde A})\Big).
\end{aligned}\end{align*}
Combining this with (\ref{sun1}), we have
\begin{align}\begin{aligned}\label{bon3}
\widetilde W_i(A,\beta)={\rm im}\Big(H_0(E^A+\beta,\;I_i(\widetilde A))\to H_0(E^A+\beta,\;S_{\widetilde A})\Big).
\end{aligned}\end{align}
Then the functoriality of $H_0(E^A+\beta,\bullet)$ induces a canonical homomorphism
$$W_i(A,\beta)={\rm im}\Big(H_0(E^A+\beta,\;I_i(  A))\to H_0(E^A+\beta,\;S_{  A})\Big)\to\widetilde W_i(A,\beta).$$
\end{proof}

\subsection{Proof of Theorem \ref{baaa}.} First, we construct the canonical epimorphism in Theorem \ref{baaa}. Second, in Proposition \ref{sun} we prove a special case of Theorem \ref{baaa}. Finally, we reduce Theorem \ref{baaa} to this special case.

\begin{proof}[Proof of part (1) in Theorem 1.4]
For any $F\in \mathfrak F_i(A)$, $I_0(F)\subset I_i(A)$. The $S_F$-module $I_0(F)$ is also an $S_A$-module via the epimorphism $S_A\to S_F$ associated to the closed immersion $\bar i_{F,\,A}\colon X_F\to X_A$. The composition $I_i(A)\hookrightarrow S_A\to\mathop{\bigoplus}\limits_{F\in\mathfrak F_i(A)}S_F$ factors as $$I_i(A)\twoheadrightarrow I_i(A)/I_{i-1}(A)\simeq\mathop{\bigoplus}\limits_{F\in\mathfrak F_i(A)}I_0(F)\hookrightarrow\mathop{\bigoplus}\limits_{F\in\mathfrak F_i(A)}S_F,$$ and it induces a commutative diagram
\begin{eqnarray}\label{431}\begin{split}\xymatrix{I_i(A)/I_{i-1}(A)\ar[r]\ar[d]^\simeq&S_A/I_{i-1}(A)\ar[d]\\
{\displaystyle \bigoplus_{F\in \mathfrak F_i(A)} I_0(F)}\ar[r]&{\displaystyle \bigoplus_{F\in \mathfrak F_i(A)} S_F}}\end{split}\end{eqnarray}
 of $S_A$-modules. Applying the right exact functor $H_0(E^A+\beta,\,\bullet)$ to the short exact sequences
\begin{eqnarray*}
&&0\to I_{i-1}(A)\to S_A\to S_A/I_{i-1}(A)\to0;\\
&&0\to I_i(A)\to S_A\to S_A/I_i(A)\to0;\\
&&0\to I_i(A)/I_{i-1}(A)\to S_A/I_{i-1}(A)\to S_A/I_i(A)\to0,
\end{eqnarray*}
we have a commutative diagram
\begin{align}\begin{aligned}\label{li1}
\xymatrix{&0\ar[d]&&H_0(E^A+\beta,\;I_i(A)/I_{i-1}(A))\ar[d]&\\
0\ar[r]&W_{i-1}(A, \beta)\ar[r]\ar[d]&H_0(E^A+\beta,\,S_A)\ar[r]\ar @{=} [d] & H_0(E^A+\beta,\;S_A/I_{i-1}(A))\ar[r]\ar[d]&0\\
0\ar[r]&W_i(A, \beta)\ar[r]&H_0(E^A+\beta,\,S_A)\ar[r]&H_0(E^A+\beta,\;S_A/I_i(A))\ar[r]&0}
\end{aligned}\end{align}
with exact rows and columns. We thus obtain
\begin{align}\begin{aligned}\label{bi1}
&\quad\;W_i(A, \beta)/W_{i-1}(A, \beta)\\&={\rm ker}\Big(H_0(E^A+\beta,\;S_A/I_{i-1}(A))\to H_0(E^A+\beta,\;S_A/I_i(A))\Big)\\
 &={\rm im}\Big(H_0(E^A+\beta,\; I_i(A)/ I_{i-1}(A))\to H_0(E^A+\beta,\; S_A/I_{i-1}(A))\Big)\\
&\twoheadrightarrow{\rm im}\Big(H_0\Big(E^A+\beta,\;\bigoplus_{F\in \mathfrak F_i(A)}I_0(F)\Big)\to H_0\Big(E^A+\beta,\;\bigoplus_{F\in \mathfrak F_i(A)}S_F\Big)\Big)\\
&=\bigoplus_{F\in \mathfrak F_i(A)}{\rm im}\Big(H_0(E^A+\beta,\;I_0(F))\to H_0(E^A+\beta,\;S_F)\Big)\\
&=\bigoplus_{F\in \mathfrak F_i(A)}{\rm im}\Big((i_{F,A})_+H_0(E^F+\beta,I_0(F))\to (i_{F, A})_+H_0(E^F+\beta,S_F)\Big)\\
&=\bigoplus_{F\in \mathfrak F_i(A)}(i_{F,\,A})_+{\rm im}\Big(H_0(E^F+\beta,\;I_0(F))\to H_0(E^F+\beta,\;S_F)\Big)\\
&=\bigoplus_{F\in \mathfrak F_i(A)}(i_{F,\,A})_+W_0(F, \beta)\\
&=\bigoplus_{\substack{F\in \mathfrak F_i(A)\\\beta\in\C F}}(i_{F,\,A})_+W_0(F, \beta).
\end{aligned}\end{align}
In this sequence of equalities, the first holds by applying the snake lemma to the last two rows of (\ref{li1}), the second uses the right column of (\ref{li1}), the third is obvious, the fourth uses Lemma \ref{1029c}, the fifth follows from the Kashiwara's equivalence for the closed immersion $i_{F,\,A}\colon \A^F\to\A^A$, the sixth follows from the definition of $W_0(F,\beta)$ and the last holds by part (1) of Lemma \ref{lao}.
Finally, the surjectivity of the third morphism follows from the commutative diagram (\ref{431}). This proves part (1) of Theorem \ref{baaa}.
\end{proof}


\begin{prop}\label{sun}
Let $F_1,\ldots,F_k$ be all facets of $A$. If for all $1\leq i\leq k$, the intersection of any $i$ distinct elements of $\{F_1,\ldots,F_k\}$ is a face of codimension $i$ in $A$, then
$$W_i(A, \beta)/W_{i-1}(A, \beta)\simeq\bigoplus_{\substack{F\in \mathfrak F_i(A)\\\beta\in\C F}}(i_{F,\,A})_+W_0(F, \beta).$$
\end{prop}

\begin{proof}
Let $F'$ be the smallest face of $A$. Then $F'=F_1\cap\cdots\cap F_k$ and $\mathbb R_{\geq0} F'=\mathbb R F'$. The condition on $A$ implies that $q(\mathbb R_{\geq0}A)$ is a simplicial cone in $\mathbb R A/\mathbb R F'$, where $q\colon \mathbb R A\to\mathbb R A/\mathbb R F'$ is the canonical epimorphism. Without loss of generality, assume that $q(\mathbb R_{\geq0}A)=\sum\limits_{j=1}^k q(\mathbb R_{\geq0} a_j)$. Let $A_0$ be the submatrix of $A$ consisting of $a_1,\ldots,a_k$ and $F'$, so that $\mathbb R_{\geq0}A=\mathbb R_{\geq0}A_0$.
We may assume that $d_A=n$. The direct sum $(\bigoplus\limits_{j=1}^k\Z a_j)\oplus \Z F'$ induces an isomorphism $\iota\colon (\bigoplus\limits_{j=1}^k\C a_j)\oplus\C F'\simeq\C^n$ such that $\iota(a_j)=e_j$ for $1\leq j\leq k$, and $\iota(F')\subset\bigoplus\limits_{j=k+1}^{n} \Z e_j$ where $e_1,\ldots,e_n$ is the canonical basis of $\C^n$. This reduces $A_0$ to be of the form
$A_0=\left(\begin{array}[c]{ll}I_k&0\\0&F'\end{array}\right).$

Consider the commutative diagram
\begin{eqnarray*}
\xymatrix{&H_0(E^A+\beta,\;I_i(A)/I_{i-1}(A))\ar[d]^\simeq\ar[r]^\varphi&H_0(E^A+\beta,\;S_A/I_{i-1}(A))\ar[d]\\
{\displaystyle\bigoplus_{F\in \mathfrak F_i(A)}}H_{1}(E^{A}+\beta,\;S_F/I_0(F))\ar[r]^{\delta}&{\displaystyle\bigoplus_{F\in \mathfrak F_i(A)}}H_0(E^A+\beta,\;I_0(F))\ar[r]^\psi\ar[ru]&{\displaystyle\bigoplus_{F\in \mathfrak F_i(A)}}H_0(E^{A}+\beta,\;S_F)}
\end{eqnarray*}
with exact rows. Recall from (\ref{bi1}) that the source and target of the epimorphism in part (1) of Theorem \ref{baaa} are isomorphic to ${\rm im}(\varphi)$ and ${\rm im}(\psi)$, respectively.
Then by a diagram chasing, to prove this epimorphism is an isomorphism, it suffices to show that for each $F\in\mathfrak F_i(A)$, the composition
\begin{eqnarray*}\label{ddaa}
H_{1}(E^{A}+\beta,\;S_F/I_0(F))\x{\delta}H_0(E^A+\beta,\;I_0(F))
\to H_0(E^{A}+\beta,\;S_A/I_{i-1}(A))
\end{eqnarray*}
is trivial. Applying Lemma \ref{4132} and Corollary \ref{ahui} to $A_0\subset A$, so we only need to show that the composition
\begin{eqnarray*}\label{daa}
H_1(E^{A_0}+\beta,\;S_F/I_0(F))\x{\delta_0} H_0(E^{A_0}+\beta,\;I_0(F))\to H_0(E^{A_0}+\beta,\;S_A/I_{i-1}(A))\end{eqnarray*}
is trivial. It is obvious if $\beta\notin\C F$, because in this case, $H_0(E^{A_0}+\beta,\;I_0(F))=0$ according to Lemma \ref{lao}. So we can assume that $\beta\in\C F=\C F_0$ where $F_0=F\cap A_0$. By the definition of the connecting map $\delta_0$, any element of ${\rm im}(\delta_0)$ is represented by an element in $\C[\p_{A_0}]\otimes_\C I_0(F)$ of the form $\sum\limits_{i=1}^n(E_i^{A_0}+\beta_i)(m_i)$ for some $m_i\in\C[\p_{A_0}]\otimes_\C S_F$. Consider the diagram
\begin{eqnarray}\label{zyou}\begin{split}\xymatrix{\C[\p_{A_0}]\otimes_\C S_F\ar[rr]^\phi\ar[d]^{E^{A_0}_i+\beta_i}&&\C[\p_{A_0}]\otimes_\C S_A/I_{i-1}(A)\ar[d]^{E^{A_0}_i+\beta_i}\\
\C[\p_{A_0}]\otimes_\C S_F\ar[rr]^\phi&&\C[\p_{A_0}]\otimes_\C S_A/I_{i-1}(A),}\end{split}\end{eqnarray}
where $\phi$ is induced by the composition $S_F\hookrightarrow S_A\twoheadrightarrow S_A/I_{i-1}(A)$. Notice that $\phi$ is $D_{F_0}$-linear but not $D_{A_0}$-linear. Obviously, $F'\prec F_0$. If ${\rm column}_i(A_0)\in F_0$, then $E_i^{A_0}+\beta_i\in D_{F_0}$ and hence the diagram (\ref{zyou}) commutes. If ${\rm column}_i(A_0)\notin F_0$, then $\beta_i=0$ and therefore $(E_i^{A_0}+\beta_i)(m_i)=0$. Consequently,
$$\phi\Big(\sum\limits_{i=1}^n(E_i^{A_0}+\beta_i)(m_i)\Big)=\sum\limits_{\substack{1\leq i\leq n\\{\rm column}_i(A_0)\in F_0}}\phi((E^{A_0}_i+\beta_i)(m_i))=\sum\limits_{\substack{1\leq i\leq n\\{\rm column}_i(A_0)\in F_0}}(E^{A_0}_i+\beta_i)(\phi(m_i)),$$
representing the zero element of $H_0(E^{A_0}+\beta,\,S_A/I_{i-1}(A))$.  This completes the proof.
\end{proof}

\begin{proof}[Proof of Theorem \ref{baaa} (2-4)]
If $A$ is normal and $\beta$ is weakly $A$-nonresonant, then $F$ is normal and $\beta$ is weakly $F$-nonresonant for any $F\prec A$ with $\beta\in\C F$. So part (3) of Theorem \ref{baaa} holds by Corollary \ref{fei}, and part (4) follows from (2) and (3). It remains to prove part (2) of Theorem  \ref{baaa}.

Recall that $F_1,\ldots,F_k$ are all facets of $A$ for which $\ell_{F_i}(\beta)\in\Z$,  and the intersection of any $i$ distinct elements of $\{F_1,\ldots,F_k\}$ is a face of $A$ of codimension $i$. Set $F_0=F_1\cap\cdots\cap F_k$ and $\epsilon_{F_0}=\sum\limits_{a_j\in F_0}a_j$.  Let ${\overline A}$ be the $n\times (N+1)$-matrix $(A,-\epsilon_{F_0})$. The map $F\mapsto\overline F:=(F,-\epsilon_{F_0})$ defines a bijection between the set of faces of $A$ containing $F_0$ and that of $\overline A$. There are two claims below.

\underline{Claim 1:} For any face $F$ of $A$, $\beta\in\C F$ implies that $F_0\prec F$.

In fact, $\beta\in\C F$ implies that $\ell_{F'}(\beta)=0$ for any facet $F'$ of $A$ containing $F$. By the assumption on $\beta$, we have $F'=F_i$ for some $1\leq i\leq k$, and then $F_0=\mathop{\bigcap}\limits_{1\leq i\leq k}F_i\prec \mathop{\bigcap}\limits_{F\prec F'\in \mathfrak F_1(A)}F'=F$. This proves Claim 1.

\underline{Claim 2:} For any $F_0\prec F\prec A$ and any $0\leq i\leq d_F$,
$H_0(E^F+\beta,\;I_i({\overline F}))=H_0(E^F+\beta,\;I_i(F))$.

In fact, $I_i({\overline F})=(I_i(F))_{t^{\epsilon_{F_0}}}$. By part (1) of Lemma \ref{lao}, we only need to show that $\beta\notin -\Z_{\geq1}\epsilon_{F_0}+{\rm qdeg}(I_i(F)/t^{\epsilon_{F_0}}I_i(F))$. If not, then $\beta=-m\epsilon_{F_0}+b+c$ for some $m\in\Z_{>0}$, $b\in\N F$ and $c\in\C^n$ such that there exists $G\prec A$ with the properties that $c\in\C G$ and $b+\N G\subseteq\deg(I_i(F)/t^{\epsilon_{F_0}}I_i(F))$. Clearly, $\epsilon_{F_0}\notin\N G$. According to $G=\mathop{\bigcap}\limits_{G\prec F'\in\mathfrak F_1(A)}F'$, there exists a facet $F'$ of $A$ such that $G\prec F'$ and $\epsilon_{F_0}\notin\N F'$. So $F'\neq F_i$ for any $1\leq i\leq k$. By the assumption on $\beta$, we have $\ell_{F'}(\beta)\notin\Z$ which contradicts to $\ell_{F'}(\beta)=\ell_{F'}(-m\epsilon_{F_0}+b+c)=\ell_{F'}(-m\epsilon_{F_0}+b)\in\Z$. This proves Claim 2.

For any face $F$ of $A$ containing $F_0$, the morphism $\pi_F\colon\A^{\overline F}\to\A^F$ induced by $F\subset\overline F$ satisfies the conditions of Corollary \ref{ahui}. We therefore have
\begin{align}\begin{aligned}\label{183201}
(\pi_F)_+W_i(\overline F,\beta)&={\rm im}\Big(H_0(E^{F}+\beta, I_i(\overline F))\to H_0(E^{F}+\beta, S_{\overline F})\Big)\\
&={\rm im}\Big(H_0(E^{F}+\beta, I_i(F))\to H_0(E^{F}+\beta, S_{F})\Big)\\
&=W_i(F, \beta),
\end{aligned}\end{align}
where the second equality holds by Claim 2. In particular, $(\pi_A)_+W_i(\overline A,\beta)=W_i(A,\beta)$. The matrix $\overline A$ satisfies the conditions of Proposition \ref{sun}, then by Proposition \ref{sun} and Claim 1,
\begin{eqnarray}\label{d1}
W_i({\overline A}, \beta)/W_{i-1}({\overline A}, \beta)&=&\bigoplus_{\substack{F\in \mathfrak F_i(A)\\F_0\prec F,\beta\in\C  F}}(i_{\overline F,\,\overline A})_+W_0(\overline F, \beta)=\bigoplus_{\substack{F\in \mathfrak F_i(A)\\\beta\in\C  F}}(i_{\overline F,\,\overline A})_+W_0(\overline F, \beta).
\end{eqnarray}
We conclude that
\begin{align*}\begin{aligned}
&\quad\;W_i(A, \beta)/W_{i-1}(A, \beta)\\
&=(\pi_A)_+W_i({\overline A},\beta)/(\pi_A)_+W_{i-1}({\overline A}, \beta)\\
&=(\pi_A)_+(W_i({\overline A},\beta)/W_{i-1}({\overline A}, \beta))\\
&=\bigoplus_{\substack{F\in \mathfrak F_i(A)\\\beta\in\C  F}}(\pi_A)_+(i_{\overline F,\,\overline A})_+W_0(\overline F, \beta)\\
&=\bigoplus_{\substack{F\in \mathfrak F_i(A)\\\beta\in\C  F}}(i_{F,\,A})_+(\pi_F)_+W_0(\overline F, \beta)\\
&=\bigoplus_{\substack{F\in \mathfrak F_i(A)\\\beta\in\C  F}}(i_{F,\,A})_+W_0(F, \beta),
\end{aligned}\end{align*}
where the first and the last equalities use (\ref{183201}), the second uses Lemma \ref{4132}, the third uses (\ref{d1}) and the fourth follows from the fact that $\pi_A\circ i_{\overline F,\,\overline A}=i_{F,A}\circ\pi_F$. This proves Theorem \ref{baaa}.
\end{proof}

\begin{rem}
(1) Claim 2 in the proof of Theorem \ref{baaa} (2-4) can be proved almost verbatim as in \cite[Corollary 5.6]{SW}. Just because the notations of this paper and \cite{SW} are slightly different, for consistency of reading, we still give a direct proof of this claim in my paper.

(2) According to part (2) of Theorem \ref{baaa} and part (2) of Lemma \ref{fei1}, if $A$ is simplicial relative to $\beta$ and $\beta\notin{\rm wRes}(A)$, we have
$$W_i(A, \beta)/W_{i-1}(A, \beta)\simeq\bigoplus_{\substack{F\in \mathfrak F_i(A)\\\beta\in\C F}}(j_{F,\,A})_{!+}\mathcal{O}_{T_F}^\beta.$$
This shows that the filtration $0\subset W_0(A,\beta)\subset\cdots\subset W_{d_A}(A,\beta)=N_A(\beta)$ on $N_A(\beta)$ has semisimple factors if $A$ is simplicial relative to $\beta$ and $\beta\notin{\rm wRes}(A)$.

It should be noted that very recently T. Reichelt and U. Walther have posted an article \cite{RW} on the arXiv where the weight filtrations for GKZ-systems are considered. They also give a filtration on $N_A(\beta)$ with semisimple factors for any normal matrix $A$ such that $\Z A=\Z^n$ and $\Z A\cap-\Z A=0$ and any $\beta\in\Z^n\backslash{\rm sRes}(A)$.

The above two filtrations coincide if $A$ and $\beta$ satisfy all the conditions mentioned above.
\end{rem}

\subsection{Counterexample to Adolphson's conjecture in the non-normal case.} Part (2) of Corollary \ref{fei} immediately implies the irreducibility conjecture of $W_0(A,\,\beta)$ in \cite{A2} for a normal matrix $A$. If $A$ is not normal, we give the following counterexample.

\begin{exam}\label{hou}
Let
\[A=\left(\begin{array}[c]{lll}1&1&0\\0&1&2\end{array}\right),\;F_2=\left(\begin{array}[c]{l}0\\2\end{array}\right).\]
Then we have a non-split short exact sequence of regular holonomic $D_A$-modules:
\begin{align}\label{ji1}0\to(j_{F_2,\,A})_{!+}\sO_{T_{F_2}}^\alpha\to W_0(A, 0)\to\widetilde  W_0(A, 0)=(j_A)_{!+}\sO_{T_A}^0\to0,\end{align}
where $\alpha=(0,1)^{\rm t}\in\C^2$. In particular, $W_0(A,0)$ is not semisimple.
\end{exam}

\begin{proof}
For the picture of $\N A$, refer to Example \ref{zy}. Let
\[\widetilde A=\left(\begin{array}[c]{llll}1&0&1&0\\0&2&1&1\end{array}\right),\;A_0=\left(\begin{array}[c]{ll}1&0\\0&2\end{array}\right),\;F_1=\left(\begin{array}[c]{l}1\\0\end{array}\right),\;\;
\widetilde F_2=\left(\begin{array}[c]{ll}0&0\\2&1\end{array}\right).\]
Then $F_1$ is a common face of $A$, $A_0$ and $\widetilde A$; $F_2$ is a common face of $A$ and $A_0$; $\widetilde F_2$ is a face of $\widetilde A$. Let $F_0$ be the common zero face of $A$, $A_0$ and $\widetilde A$. The inclusions $A_0\subset A\subset\widetilde A$ induce two projections
$$\A^{\widetilde A}={\rm Spec}\,\C[x_1,x_2,x_3,x_4]\x{\pi}\A^A={\rm Spec}\,\C[x_1,x_2,x_3]\x{p}\A^{A_0}={\rm Spec}\,\C[x_1,x_2].$$

We prove this example in 5 steps. In Step (i), we compute the composition factors of $N_A(0)$ and $(j_A)_+\sO_{T_A}^0$. In Step (ii), we prove that ${\rm ker}\Big(W_0(A, \,0)\to\widetilde W_0(A, \,0)\Big)=H_1(E^A,\;S_{\widetilde A}/S_A)$. In Step (iii), we prove that $H_1(E^A,\;S_{\widetilde A}/S_A)=(j_{F_2,\,A})_{!+}\sO_{T_{F_2}}^\alpha$. In Step (iv), we prove that the non-semisimplicity of $W_0(A,0)$ follows from that of $p_+W_0(A,0)$. In Step (v), we prove that $p_+W_0(A,0)$ is not semisimple.

\underline{Step (i).} In this step, we compute the composition factors on $N_A(0)$ and $(j_A)_+\sO_{T_A}^0$.

Recall that for any face $F$ of $A$, $j_{F,A}:T_F\to\A^A$ is the morphism induced by the homomorphism
$$\C[x_A]\to\C[t^{\pm a_1},\ldots,t^{\pm a_N}],\;\;x_j\mapsto t^{a_j}\hbox{ for any }1\leq j\leq N.$$
Every proper face of $A$ or $\widetilde A$ is normal, and $0$ is weakly nonresonant with respect to each face of $A$ or $\widetilde A$. In other words, $A$ and $\widetilde A$ satisfy all conditions in Theorem \ref{baaa} (4) for $i=1,2$.
So we have
\begin{eqnarray}\label{1831711}
&&W_1(A,0)/W_0(A, 0)=(j_{ F_1,\,A})_{!+}\sO_{T_{F_1}}^0\oplus(j_{ F_2,\,A})_{!+}\sO_{T_{F_2}}^0;\\\label{183174}
&&W_2(A,0)/W_1(A,0)=(j_{F_0,A})_+\C;\\\label{183172}
&&W_1(\widetilde A,0)/W_0(\widetilde A,0)=(j_{F_1,\,\widetilde A})_{!+}\sO_{T_{F_1}}^0\oplus(j_{\widetilde F_2,\,\widetilde A})_{!+}\sO_{T_{\widetilde F_2}}^0;\\\label{183175}
&&W_2(\widetilde A,0)/W_1(\widetilde A,0)=(j_{F_0,\widetilde A})_+\C.
\end{eqnarray}
Then
\begin{align}\begin{aligned}\label{sun3}
&\;\;\widetilde W_1(A,0)/\widetilde W_0(A, 0)\\
=&\;\;\pi_+\big(W_1(\widetilde A,0)/W_{0}(\widetilde A,0)\big)\\
=&\;\;\pi_+(j_{F_1,\,\widetilde A})_{!+}\sO_{T_{F_1}}^0\oplus\pi_+(j_{\widetilde F_2,\,\widetilde A})_{!+}\sO_{T_{\widetilde F_2}}^0\\
=&\;\;(j_{F_1,\,A})_{!+}\sO_{T_{F_1}}^0\oplus(j_{F_2,\,A})_{!+}(\pi_{F_2,\,\widetilde F_2})_+\sO_{T_{\widetilde F_2}}^0\\
=&\;\;(j_{F_1,\,A})_{!+}\sO_{T_{F_1}}^0\oplus(j_{F_2,\,A})_{!+}\sO_{T_{F_2}}^0\oplus(j_{F_2,\,A})_{!+}\sO_{T_{F_2}}^\alpha,
\end{aligned}\end{align}
where the first equality uses Lemma \ref{18317}, the second uses (\ref{183172}), the third holds by applying Lemma \ref{4131} to the following commutative diagram
\[\xymatrix{T_{F_1}\ar[r]^{j_{F_1,\widetilde A}}\ar[dr]_{j_{F_1,A}}&\A^{\widetilde A}\ar[d]^{\pi}&T_{\widetilde F_2}\ar[l]_{j_{\widetilde F_2,\,\widetilde A}}\ar[d]^{\pi_{F_2,\,\widetilde F_2}}\\
&\A^A&T_{F_2}\ar[l]_{j_{F_2,A}},}\]
and the last follows from the fact that $(\pi_{F_2,\,\widetilde F_2})_+\sO_{T_{\widetilde F_2}}^0=\sO_{T_{F_2}}^0\oplus\sO_{T_{F_2}}^\alpha.$
Similarly, we also have
\begin{eqnarray}\label{li2}
\widetilde W_2(A,0)/\widetilde W_1(A,0)=\pi_+(W_2(\widetilde A,0)/W_1(\widetilde A,0))=\pi_+(j_{F_0,\widetilde A})_+\C=(j_{F_0,A})_+\C.
\end{eqnarray}

\underline{Step (ii).} ${\rm ker}\Big(W_0(A, \,0)\to\widetilde W_0(A, \,0)\Big)=H_1(E^A,\;S_{\widetilde A}/S_A)$.

Recall from Lemma \ref{18317} that there is a canonical homomorphism $W_i(A,0)\to\widetilde W_i(A,0)$. Then (\ref{1831711}) and (\ref{sun3}) give a natural monomorphism $W_1(A,0)/W_{0}(A, 0)\to\widetilde W_1(A,0)/\widetilde W_{0}(A,0)$, which follows that
\begin{align}\label{bon1}{\rm ker}\Big(W_0(A, 0)\to\widetilde W_0(A, 0)\Big)={\rm ker}\Big(W_1(A, 0)\to\widetilde W_1(A, 0)\Big).\end{align}
By (\ref{183174}) and (\ref{li2}), $ W_2(A,0)/W_1(A,0)=\widetilde W_2(A,0)/\widetilde W_1(A,0),$
which immediately implies that
\begin{align}\label{li3}{\rm ker}\Big(W_1(A, 0)\to\widetilde W_1(A, 0)\Big)={\rm ker}\Big(W_2(A,\,0)\to\widetilde W_2(A,\,0)\Big),\end{align}
Taking $i=2$ in (\ref{bon3}) of Lemma \ref{18317}, we get $\widetilde W_2(A,0)=H_0(E^A,S_{\widetilde A}).$ This, together with (\ref{bon1}) and (\ref{li3}), implies that
\begin{align}\label{bo2}
{\rm ker}\Big(W_0(A, \,0)\to\widetilde W_0(A, \,0)\Big)={\rm ker}\Big(H_0(E^A,S_A)\to H_0(E^A,S_{\widetilde A})\Big).
\end{align}
Applying $\sK_\bullet(E^A,\bullet)$ to the short exact sequence of toric $\C[x_A]$-modules:
$$0\to S_A\to S_{\widetilde A}\to S_{\widetilde A}/S_A\to0,$$
we get the following exact sequence in ${\rm Mod}_{rh}^{X_A}(D_A)$:
\begin{eqnarray}\label{183204}
H_1(E^A,\;S_{\widetilde A})\to H_1(E^A,\;S_{\widetilde A}/S_A)\to H_0(E^A,\;S_A)\to H_0(E^A,\;S_{\widetilde A}).
\end{eqnarray}
By Corollary \ref{ahui}, $H_1(E^A,\;S_{\widetilde A})=\pi_+H_1(E^{\widetilde A},\;S_{\widetilde A})$, and by \cite[ Corollary 3.8]{SW}, $H_i(E^{\widetilde A},\;S_{\widetilde A})=0$ for any $i>0$. Consequently, the first term
of the exact sequence (\ref{183204}) vanishes, and so
\begin{align*}
{\rm ker}\Big(H_0(E^A,S_A)\to H_0(E^A,S_{\widetilde A})\Big)=H_1(E^A,\;S_{\widetilde A}/S_A).
\end{align*}
Combining this with (\ref{bo2}), we prove Step (ii).

\underline{Step (iii).} $H_1(E^A,\;S_{\widetilde A}/S_A)=(j_{F_2,\,A})_{!+}\sO_{T_{F_2}}^\alpha$ and (\ref{ji1}) is exact.

First note that the $S_A$-module $S_{\widetilde A}/S_A=\C[t_1,t_2]/\C[t_1,t_1t_2,t_2^2]$ has a natural toric $S_{F_2}=\C[t_2^2]$-module structure, which is isomorphic to $t_2\C[t_2^2]\simeq S_{F_2}[-\alpha]$.
So, by Lemma \ref{1029c},
\begin{eqnarray}\label{12122}H_1(E^A,\;S_{\widetilde A}/S_A)=(i_{F_2,A})_+ H_1(E^{F_2},\;S_{F_2}[-\alpha])=(i_{F_2,A})_+ H_1(E^{F_2}-\alpha,\;S_{F_2}).\end{eqnarray}
Let $F_2'$ be the $1\times 1$-matrix $(2)$. By the proof of the equality (\ref{tian1}), we have
\begin{align}\begin{aligned}\label{b5}
H_1(E^{F_2}-\alpha,\;S_{F_2})&=H_0(E^{F'_2}-1,S_{F'_2})\oplus H_1(E^{F'_2}-1,S_{F'_2})\\
&=H_0(E^{F_2}-1,S_{F_2})\oplus H_1(E^{F'_2}-1,S_{F'_2}).
\end{aligned}\end{align}
Recall that $j_{F_2}:T_{F_2}={\rm Spec}\,\C[t_2^{\pm2}]\to\A^{F_2}={\rm Spec}\,\C[x_2]$ is the morphism induced by the homomorphism $\C[x_2]\to\C[t_2^{\pm2}],\;x_2\mapsto t_2^2$. As $1\notin{\rm Res}(F_2')=2\Z$, then by \cite[Corollary 3.8]{SW}, $H_1(E^{F'_2}-1,S_{F'_2})=0$, and by part (3) of Lemma \ref{fei1},
\begin{align}\label{r1}H_0(E^{F_2}-\alpha,S_{F_2})=(j_{F_2})_{!+}\sO_{T_{F_2}}^{-\alpha}=(j_{F_2})_{!+}\sO_{T_{F_2}}^\alpha.\end{align}
\if<\begin{align}\begin{aligned}\label{bo7}
&H_0\Big(\C[\p_{2}]\otimes_\C \C[t_2^2]\x{2\p_{2}x_2-1}\C[\p_{2}]\otimes_\C \C[t_2^2]\Big)=(j_{F_2})_{!+}\sO_{T_{F_2}}^\alpha;\\
&H_1\Big(\C[\p_{2}]\otimes_\C \C[t_2^2]\x{2\p_{2}x_2-1}\C[\p_{2}]\otimes_\C \C[t_2^2]\Big)=0.
\end{aligned}\end{align}
As $E_1^{F_2}=0$ and $E_2^{F_2}=2\p_{2}x_2$, we have
\begin{align*}\begin{aligned}
&\quad H_1(E^{F_2},\;S_{\widetilde A}/S_A)\\
&=H_1\sK\Big(E_1^{F_2},E_2^{F_2};\C[\p_2]\otimes_\C t_2\C[t_2^2]\Big)\\
&=\bigoplus_{i=0}^1H_i\Big(\C[\p_{2}]\otimes_\C t_2\C[t_2^2]\x{2\p_{2}x_2}\C[\p_{2}]\otimes_\C t_2\C[t_2^2]\Big)\\
&=\bigoplus_{i=0}^1H_i\Big(\C[\p_{2}]\otimes_\C \C[t_2^2]\x{2\p_{2}x_2-1}\C[\p_{2}]\otimes_\C \C[t_2^2]\Big)\\&=(j_{F_2})_{!+}\sO_{T_{F_2}}^\alpha.&&\text{by (\ref{bo7})}
\end{aligned}\end{align*}>\fi
Substituting these into (\ref{12122}) and (\ref{b5}), we have
\begin{eqnarray}\label{12123}H_1(E^A,\,S_{\widetilde A}/S_A)=(i_{F_2,\,A})_+(j_{F_2})_{!+}\sO_{T_{F_2}}^\alpha=(j_{F_2,\,A})_{!+}\sO_{T_{F_2}}^\alpha.\end{eqnarray}
Then the exactness of (\ref{ji1}) follows immediately from Step (ii), (\ref{12123}) and the irreducibility of $\widetilde W_0(A, 0)=(j_A)_{!+}\sO_{T_A}^0$. This proves Step (iii).

\underline{Step (iv).} If $W_0(A,0)$ is a semisimple $D_A$-module, then $p_+W_0(A,0)$ is a semisimple $D_{A_0}$-module.

To prove it, suppose that $W_0(A,0)$ is semisimple. Then (\ref{ji1}) is a split exact sequence, and so
\begin{align*}
p_+W_0(A,0)&=p_+(j_{F_2,\,A})_{!+}\sO_{T_{F_2}}^\alpha\oplus p_+ (j_A)_{!+}\sO_{T_A}^0\\
&=(j_{F_2,A_0})_{!+}\sO_{T_{F_2}}^\alpha\oplus(j_{A_0})_{!+}(\pi_{A_0,A})_+\sO_{T_A}^0\\
&=(j_{F_2,A_0})_{!+}\sO_{T_{F_2}}^\alpha\oplus(j_{A_0})_{!+}\sO_{T_{A_0}}^0\oplus(j_{A_0})_{!+}\sO_{T_{A_0}}^\alpha,
\end{align*}
where the second equality holds for the same reason as in the third equality in (\ref{sun3}), and the third follows because the homomorphism $\pi_{A_0,A}:T_A\to T_{A_0}$ of tori induces that $(\pi_{A_0,A})_+\sO_{T_A}^0=\sO_{T_{A_0}}^0\oplus\sO_{T_{A_0}}^\alpha$.
This proves Step (iv).

\underline{Step (v).} The $D_{A_0}$-module $p_+W_0(A,0)$ is not semisimple.

First note that there are two isomorphisms
\begin{eqnarray}\begin{split}\label{12223}
I_0(A)&=t_1t_2\C[t_1,t_2]\simeq t_1t_2\C[t_1,t_2^2]\oplus t_1t_2^2\C[t_1,t_2^2];\\
S_A&=\C[t_1,t_1t_2,t_2^2]=t_1t_2\C[t_1,t_2^2]\oplus\C[t_1,t_2^2]
\end{split}\end{eqnarray}
of toric $\C[x_{A_0}]=\C[x_1,x_2]\simeq\C[t_1,t_2^2]$-modules, which are compatible with the inclusion $I_0(A)\subset S_A$. For any $D_{F_1}=\C[x_1,\p_1]$-module $M_1$ and any $D_{F_2}=\C[x_2,\p_2]$-module $M_2$, denote by $M_1\boxtimes M_2$ the induced $D_{A_0}=\C[x_1,x_2,\p_1,\p_2]$-module $M_1\otimes_\C M_2$. Let $\alpha'=(1,0)^{\rm t}$. We have
\begin{align*}
&\quad\; p_+W_0(A,0)\\
&={\rm im}\Big(H_0(E^{A_0},I_0(A))\to H_0(E^{A_0},S_A)\Big)\\
&={\rm im}\Big(H_0(E^{A_0},t_1t_2\C[t_1,t_2^2]\oplus t_1t_2^2\C[t_1,t_2^2])\to H_0(E^{A_0},t_1t_2\C[t_1,t_2^2]\oplus\C[t_1,t_2^2])\Big)\\
&=H_0(E^{A_0},t_1t_2\C[t_1,t_2^2])\oplus{\rm im}\Big(H_0(E^{A_0},t_1t_2^2\C[t_1,t_2^2])\to H_0(E^{A_0},\C[t_1,t_2^2])\Big)\\
&=\Big(H_0(E^{F_1},t_1\C[t_1])\boxtimes H_0(E^{F_2},t_2\C[t_2^2])\Big)\oplus\\&\quad\;\Big({\rm im}\Big(H_0(E^{F_1},t_1\C[t_1])\to H_0(E^{F_1},\C[t_1])\Big)\boxtimes {\rm im}\Big(H_0(E^{F_2},t_2^2\C[t_2^2])\to H_0(E^{F_2},\C[t_2^2])\Big)\Big)\\
&=\Big(H_0(E^{F_1}-\alpha',\C[t_1])\boxtimes H_0(E^{F_2}-\alpha,\C[t_2^2])\Big)\oplus\Big(W_0(F_1,0)\boxtimes W_0(F_2,0)\Big)\\
&=\Big(D_{F_1}/D_{F_1}(\p_1x_1-1)\boxtimes (j_{F_2})_{!+}\sO_{T_{F_2}}^\alpha\Big)\oplus\Big((j_{F_1})_{!+}\sO_{T_{F_1}}^0\boxtimes (j_{F_2})_{!+}\sO_{T_{F_2}}^0\Big).
\end{align*}
In this sequence of equalities, the first uses Corollary \ref{ahui}, the second uses (\ref{12223}), the third and fourth are obvious, the fifth follows from (2) of Definition \ref{dimao} and the last follows from (\ref{r1}) and part (2) of Corollary \ref{fei}. So, to prove the non-semisimplicity for $p_+W_0(A,0)$, we only need to show that for $D_{F_1}/D_{F_1}(\p_1x_1-1)=D_{F_1}/D_{F_1}x_1\p_1$ or its Fourier transform $D_{F_1}/D_{F_1}\p_1x_1$.
Note that $j_{F_1}\colon T_{F_1}={\rm Spec}\,\C[x_1,x_1^{-1}]\to\A^{F_1}={\rm Spec}\,\C[x_1]$ is an open immersion with complement $j_{F_0,F_1}\colon{\rm Spec}\,\C\to{\rm Spec}\,\C[x_1]$. Applying Theorem \ref{baaa} to $F_1$, we have a short exact sequence
\begin{eqnarray*}\label{18314}0\to (j_{F_1})_{!+}\sO_{T_{F_1}}^0\to(j_{F_1})_+\sO_{T_{F_1}}^0=D_{F_1}/D_{F_1}\p_1x_1\to (j_{F_0,F_1})_{+}\C\to0,
\end{eqnarray*}
which is non-split according to
$${\rm Hom}_{D_{F_1}}((j_{F_0,F_1})_{+}\C,(j_{F_1})_+\sO_{T_{F_1}}^0)={\rm Hom}_{D_{T_{F_1}}}(j_{F_1}^+(j_{F_0,F_1})_{+}\C,\sO_{T_{F_1}}^0)=0.$$
This completes the proof.
\end{proof}

\section{Proof of Theorem \ref{yaode}}
In this section, we prove Theorem \ref{yaode}, a perverse sheaf analog of Theorem \ref{baaa}. In $\S 5.1$, we prove some functorial and vanishing properties for perverse sheaves on $X_A$. In $\S 5.2$, we prove Theorem \ref{yaode} in five steps. In $\S 5.3$, we give a counterexample of Theorem \ref{yaode} without the simplicial condition.

Recall that the open immersion $\bar j_A:T_A\to X_A$ induced by the inclusion $\C[t^{a_1},\ldots,t^{a_N}]\subset\C[t^{\pm a_1},\ldots,t^{\pm a_N}]$ factors as
$T_A\x{\bar \ell_i^A}U_i(A)=X_A-\mathop{\bigcup}\limits_{\substack{F\prec A\\d_F<d_A-i}}X_F\x{\bar k_i^A}X_A.$ Let $\sL_A$ is a rank one local system on $T_A$. Then $\sW_i(\sL_A)=(\bar k_i^A)_{!*}(\bar \ell_i^A)_*(\sL_A[d_A])$ is a perverse subsheaf of $(\bar j_A)_*(\sL_A[d_A])$.

Denote by $\pi_{B,C}$ the induced homomorphism $T_C\to T_B$ of tori for any integer matrix $C$ and any subset $B$ of $C$. Recall that for any local system $\sL$ on $T_C$, $\pi_{B,C}^{\diamond}(\sL )$ denotes the set of isomorphism classes of local systems on $T_B$ whose inverse images on $T_C$ are isomorphic to $\sL$. For any set $S$ of isomorphism classes of local systems on $T_C$, denote $\pi_{C,B}^\diamond(S)$ in the same sense.

\subsection{Functorial and vanishing properties of perverse sheaves on $X_A$.}
Lemma \ref{4134} ensures that we can reduce Theorem \ref{yaode} to the case when $A$ is normal. Lemma \ref{1214} calculates the restriction of $\sW_i(\sL_A)$ to orbits of the toric variety $X_A$. Lemma \ref{1215} will be needed to construct the canonical epimorphism in Theorem \ref{yaode}.
\begin{lem}\label{4134}
Let $B$ be a subset of $A$ such that $\mathbb R_{\geq 0}A=\mathbb R_{\geq0}B$. Then $B\subset A$ induces a commutative diagram of schemes:
\[\xymatrix{T_A\ar[r]^{\bar \ell_i^A}\ar[d]^{\pi_{B,\,A}}&U_i(A)\ar[r]^{\bar k_i^A}\ar[d]^{ \pi_i}&X_A\ar[d]^{\pi}\\T_{B}\ar[r]^{\bar \ell_i^{B}}&U_i(B)\ar[r]^{\bar k_i^{B}}&X_{B}.}\]
Moreover,
\begin{eqnarray*}
(\pi_{B,\,A})_*(\sL_A)&=&\bigoplus_{\substack{\sL_B\in\pi_{B,A}^{\diamond}(\sL_A)}}\sL_B;\\
\pi_*(\sW_i(\sL_A)/\sW_{i-1}(\sL_A))&=&\bigoplus_{\substack{\sL_B\in\pi_{B,A}^{\diamond}(\sL_A)}}(\sW_i(\sL_B)/\sW_{i-1}(\sL_B)).
\end{eqnarray*}
\end{lem}

\begin{proof}
By $\mathbb R_{\geq0}A=\mathbb R_{\geq0}B$, the group $\Z A/\Z B$ is finite. For some isomorphisms $ \Z A\simeq\Z^{d_A}$ and $\Z^{d_A}\simeq\Z B$, the composition $\Z^{d_A}\simeq\Z B\subset\Z A\simeq\Z^{d_A}$ is defined by a diagonal matrix ${\rm diag}\{d_1,\,\ldots,d_{d_A}\}$. So we can identify $T_A\simeq \G_{\rm m}^{d_A}$ and $\G_{\rm m}^{d_A}\simeq T_B$ such that $\pi_{B,\,A}\colon T_A\to T_B$ corresponds to $\prod\limits_{i=1}^{d_A}p_i\colon  \G_{\rm m}^{d_A}\to\G_{\rm m}^{d_A}$ where $p_i\colon \G_{\rm m}\to \G_{\rm m}$ is the $d_i$-th power map. We have $\sL_A\simeq\sL_1\boxtimes\cdots\boxtimes\sL_{d_A}$ for some rank one local systems $\sL_1,\ldots,\sL_{d_A}$ on $\G_{\rm m}$. The local system $\sL_i$ is defined by a multi-valued function $z^{\alpha_i}$ on $\C^*$ for some $\alpha_i\in\C$. Write $\sL_i$ as $[z^{\alpha_i}]$. Then
$$(p_i)_*\sL_i=(p_i)_*[z^{\alpha_i}]=\bigoplus_{e_i=0}^{d_i-1}[z^{\frac{\alpha_i+e_i}{d_i}}]$$
and hence
$$(\pi_{B,\,A})_*\sL_A=\mathop{\boxtimes}\limits_{i=1}^{d_A}(p_i)_*\sL_i=\mathop{\boxtimes}\limits_{i=1}^{d_A}\Big(\bigoplus_{e_i=0}^{d_i-1}[z^{\frac{\alpha_i+e_i}{d_i}}]\Big)=\bigoplus_{\substack{0\leq e_i<d_i}}\Big(\mathop{\boxtimes}\limits_{i=1}^{d_A}[z^{\frac{\alpha_i+e_i}{d_i}}]\Big)=\bigoplus_{\substack{\sL_B\in\pi_{B,A}^{\diamond}(\sL_A)}}\sL_B.$$
According to Lemma \ref{4132}, $\pi$ and $\pi_i$ are finite morphisms. Then we have
\begin{align*}
\pi_*\sW_i(\sL_A)&=\pi_*(\bar k_i^A)_{!*}(\bar \ell_i^A)_*(\sL_A[d_A])\\
&=(\bar k_i^{B})_{!*}(\pi_i)_*(\bar \ell_i^A)_*(\sL_A[d_A])\\
&=(\bar k_i^{B})_{!*}(\bar \ell_i^{B})_*(\pi_{B,\,A})_*(\sL_A[d_A])\\
&=\bigoplus_{\substack{\sL_B\in\pi_{B,A}^{\diamond}(\sL_A)}}(\bar k_i^{B})_{!*}(\bar \ell_i^{B})_*(\sL_{B}[d_{B}])\\
&=\bigoplus_{\substack{\sL_B\in\pi_{B,A}^{\diamond}(\sL_A)}}\sW_i(\sL_{B}),
\end{align*}
where the second equality follows from the proof of part (2) in Lemma \ref{4131}. So the exact functor $\pi_*:{\rm Perv}(X_A)\to{\rm Perv}(X_B)$ immediately implies the second equation of this lemma.
\end{proof}

Recall that for any $F\prec A$, $\bar i_{F,A}:X_F\to X_A$ is the closed immersion defined by the ideal of $S_A$ generated by those $t^{a_j}$ such that $a_j\notin F$, and $\bar j_{F,A}=\bar i_{F,A}\circ \bar j_F$ is the orbit embedding of $T_F$ into the toric variety $X_A$. We calculate the restrictions of $\sW_i(\sL_{A})$ on orbits of $X_A$ as follows.
In this paper, we follow the convention that $\textstyle\bigwedge^\bullet\C^k$ lives in cohomological degrees $0$ through $k$.
\begin{lem}\label{1214}
Suppose that  $A$ is normal. Let $F$ be a face of $A$, and let $\sL_A$ be a rank one local system on $T_A$. If there is a rank one local system $\sL_F$ on $T_F$ such that $\pi_{F,A}^{-1}(\sL_F)=\sL_A$, then
\begin{eqnarray}\label{ji4}
\label{bi}&&\bar j_{F,\,A}^{-1}(\bar j_A)_*(\sL_A[d_A])=\sL_F[d_F]\otimes_\C\textstyle\bigwedge^\bullet\C^{d_A-d_F}[d_A-d_F]\in D_c^b(T_F).
\end{eqnarray}

Otherwise,  for any $0\leq i\leq d_A$ we have
\begin{eqnarray}\label{zhen}
\bar j_{F,\,A}^{-1}\sW_i(\sL_{A})=\bar j_{F,\,A}^!\sW_i(\sL_{A})=0.
\end{eqnarray}
\end{lem}

\begin{proof}
Since $A$ is normal, the short exact sequence
$$0\to\Z F\to\Z A\x{\tau}\Z A/\Z F\to0$$
splits. Fix a section $\iota$ of $\tau$ and set $\epsilon_F=\sum\limits_{a_j\in F}a_j$. By the map
$$G\mapsto\overline G:=(G,-\epsilon_F)\mapsto\underline G:=\iota(\tau(G)),$$ the three sets \{faces of $A$ containing $F$\}, \{faces of $\overline A:=(A,-\epsilon_F)$\} and \{faces of $\underline A:=\iota(\tau(A))$\} are in bijective correspondence. The section $\iota$ also defines a commutative diagram
\begin{eqnarray*}
\xymatrix{
T_F\times_\C T_{\underline A}\ar[rr]^{{\rm id}_{T_F}\times \bar \ell_i^{\underline A}}\ar[d]^\simeq&&T_F\times_\C U_i(\underline A)\ar[rr]^{{\rm id}_{T_F}\times \bar k_i^{\underline A}}\ar[d]^\simeq&&T_F\times_\C X_{\underline A}\ar[d]^\simeq&&T_F\times_\C T_{\underline F}\ar[ll]_{{\rm id}_{T_F}\times \bar j_{\underline F,\underline A}}\ar[d]^\simeq\\
T_{\overline A}\ar[rr]^{\bar \ell_i^{\overline A}}\ar @{=} [d] &&U_i(\overline A)\ar[rr]^{\bar k_i^{\overline A}}\ar[d]^{j_i}&&X_{\overline A}\ar[d]^j&&T_{\overline F}\ar @{=} [d] \ar[ll]_{\bar j_{\overline F,\overline A}}\\
T_{  A}\ar[rr]^{\bar \ell_i^{  A}}&&U_i(A)\ar[rr]^{\bar k_i^{ A}}&&X_{ A}&&T_F\ar[ll]_{\bar j_{ F,  A}}.}\end{eqnarray*}
Here $j$ and $j_i$ are the open immersions induced by $A\subset\overline A$, and $\bar j_{\underline F,\,\underline A}\colon T_{\underline F}={\rm Spec}\,\C\to X_{\underline A}$ is defined by the unique $T_{\underline A}$-invariant point of $X_{\underline A}(\C)$.

As $T_A=T_{\overline A}$, $\sL_A$ can be viewed as a local system $\sL_{\overline A}$ on $T_{\overline A}=T_F\times T_{\underline A}$. There exists a unique local system $\sL_F$ on $T_F$ and a unique local system $\sL_{\underline A}$ on $T_{\underline A}$ such that $\sL_{\overline A}=\sL_F\boxtimes\sL_{\underline A}$.
By the above commutative diagram, we have
\begin{align}\begin{aligned}\label{wen}
\sW_i(\sL_{\overline A})&=(\bar k_i^{\overline A})_{!*}(\bar \ell_i^{\overline A})_*(\sL_{\overline A}[d_{\overline A}])\\
&=({\rm id}_{T_F}\times\bar k_i^{\underline A})_{!*}({\rm id}_{T_F}\times\bar \ell_i^{\underline A})_*(\sL_F[d_F]\boxtimes\sL_{\underline A}[d_{\underline A}])\\
 &=\sL_F[d_F]\boxtimes(\bar k_i^{\underline A})_{!*}(\bar \ell_i^{\underline A})_*(\sL_{\underline A}[d_{\underline A}])\\
&=\sL_F[d_F]\boxtimes\sW_i(\sL_{\underline A}),
\end{aligned}\end{align}
and therefore
\begin{align}\begin{aligned}\label{qiao}
\bar j_{F,\,A}^{-1}\sW_i(\sL_A)&=\bar j^{-1}_{\overline F,\,\overline A}j^{-1}(\bar k_i^A)_{!*}(\bar \ell_i^A)_*(\sL_A[d_A])\\
&=\bar j^{-1}_{\overline F,\,\overline A}(\bar k_i^{\overline A})_{!*}j_i^{-1}(\bar\ell_i^A)_*(\sL_A[d_A])\\
&=\bar j^{-1}_{\overline F,\,\overline A}(\bar k_i^{\overline A})_{!*}(\bar \ell_i^{\overline A})_*(\sL_{\overline A}[d_{\overline A}])\\
&=\bar j^{-1}_{\overline F,\,\overline A}\big(\sL_F[d_F]\boxtimes\sW_i(\sL_{\underline A})\big)\\
&=\sL_F[d_F]\boxtimes\bar j_{\underline F,\,\underline A}^{-1}\sW_i(\sL_{\underline A}).
\end{aligned}\end{align}
The action $\rho$ of the torus $T_{\underline A}$ on $X_{\underline A}$ defines the following cartesian diagram
\[\xymatrix{ T_{\underline A}\times T_{\underline A}\ar[rr]^{{\rm id}_{T_{\underline A}}\times\bar \ell_i^{\underline A}}\ar[d]^m&&T_{\underline A}\times U_i(\underline A)\ar[rr]^{{\rm id}_{T_{\underline A}}\times\bar k_i^{\underline A}}\ar[d]^{\rho_i}&&T_{\underline A}\times X_{\underline A}\ar[d]^{\rho}\\
 T_{\underline A}\ar[rr]^{\bar \ell_i^{\underline A}}&&U_i(\underline A)\ar[rr]^{\bar k_i^{\underline A}}&& X_{\underline A},}\]
where $m$ is the multiplicative map and $\rho_i$ is the restriction of $\rho$. Then $\rho$, $\rho_i$ and $m$ are smooth morphisms of relative dimension $d_{\underline A}$, and $m^{-1}\sL_{\underline A}=\sL_{\underline A}\boxtimes\sL_{\underline A}$. Hence $\rho^{-1}[d_{\underline A}]=\rho^![-d_{\underline A}]$, $\rho_i^{-1}[d_{\underline A}]=\rho_i^![-d_{\underline A}]$ and $m^{-1}[d_{\underline A}]=m^![-d_{\underline A}]$ are exact functors of perverse sheaves. Consequently,
\begin{align}\begin{aligned}\label{18322}
(\rho^{-1}[d_{\underline A}])\sW_i(\sL_{\underline A})
&=(\rho^{-1}[d_{\underline A}])(\bar k_i^{\underline A})_{!*}(\bar \ell_i^{\underline A})_*(\sL_{\underline A}[d_{\underline A}])\\
&=({\rm id}_{T_{\underline A}}\times\bar k_i^{\underline A})_{!*}(\rho_i^{-1}[d_{\underline A}])(\bar \ell_i^{\underline A})_*(\sL_{\underline A}[d_{\underline A}])\\
&=({\rm id}_{T_{\underline A}}\times\bar k_i^{\underline A})_{!*}({\rm id}_{T_{\underline A}}\times\bar \ell_i^{\underline A})_*(m^{-1}[d_{\underline A}])(\sL_{\underline A}[d_{\underline A}])\\
&=({\rm id}_{T_{\underline A}}\times\bar k_i^{\underline A})_{!*}({\rm id}_{T_{\underline A}}\times\bar \ell_i^{\underline A})_*(\sL_{\underline A}[d_{\underline A}]\boxtimes\sL_{\underline A}[d_{\underline A}])\\
&=\sL_{\underline A}[d_{\underline A}]\boxtimes(\bar k_i^{\underline A})_{!*}(\bar \ell_i^{\underline A})_*(\sL_{\underline A}[d_{\underline A}])\\
&=\sL_{\underline A}[d_{\underline A}]\boxtimes\sW_i(\sL_{\underline A}),
\end{aligned}\end{align}
where the second equality follows from the exactness of $\rho^{-1}[d_{\underline A}]$ and $\rho_i^{-1}[d_{\underline A}]$, and the third follows from that of $\rho_i^{-1}[d_{\underline A}]$ and $m^{-1}[d_{\underline A}]$. Considering the commutative diagram
\[\xymatrix{T_{\underline A}=T_{\underline A}\times{\rm Spec}\,\C\ar[d]^{pr}\ar[rr]^{{\rm id}_{T_{\underline A}}\times \bar j_{\underline F,\underline A}}&&T_{\underline A}\times T_{\underline A}\ar[d]^\rho\\
{\rm Spec}\,\C\ar[rr]^{\bar j_{\underline F,\underline A}}&&X_{\underline A},}\]
where $pr\colon T_{\underline A}\to{\rm Spec}\,\C$ is the structure morphism, then by (\ref{18322}) we have
\begin{align}\begin{aligned}\label{183142}
pr^{-1}\bar j_{\underline F,\,\underline A}^{-1}\sW_i(\sL_{\underline A})&=({\rm id}_{T_{\underline A}}\times \bar j_{\underline F,\,\underline A})^{-1}(\rho^{-1}[d_{\underline A}])(\sW_i(\sL_{\underline A}))[-d_{\underline A}]\\
&=({\rm id}_{T_{\underline A}}\times \bar j_{\underline F,\,\underline A})^{-1}(\sL_{\underline A}[d_{\underline A}]\boxtimes\sW_i(\sL_{\underline A}))[-d_{\underline A}]\\
&=\sL_{\underline A}\boxtimes \bar j_{\underline F,\,\underline A}^{-1}\sW_i(\sL_{\underline A}).
\end{aligned}\end{align}

If $\pi_{F,A}^{-1}(\sL_F)\simeq\sL_A$, then $\sL_{\underline A}$ is the constant sheaf $\C_{T_{\underline A}^{\rm an}}$. For any sufficient small analytic open neighborhood $U$ of $0$ in $X_{\underline A}^{\rm an}$, $U\cap T_{\underline A}^{\rm an}$ is homotopic to $(\C^*)^{d_{\underline A}}$.
For such $U$, we have
\begin{eqnarray}\label{diao}\bar j_{\underline F,\,\underline A}^{-1}(\bar j_{\underline A})_*\sL_{\underline A}=\bar j_{\underline F,\,\underline A}^{-1}(\bar j_{\underline A})_*\C_{T_{\underline A}^{\rm an}}=R\Gamma(U\cap T_{\underline A}^{\rm an},\;\C)=R\Gamma((\C^*)^{d_{\underline A}},\;\C)=\textstyle\bigwedge^\bullet\C^{d_A-d_F}.\end{eqnarray}
According to (\ref{qiao}) and (\ref{diao}), we get (\ref{bi}).

If $\pi_{F,A}^{-1}(\sL_F)\not\simeq\sL_A$, then $\sL_{\underline A}$ is not a constant sheaf. The left hand side of (\ref{183142}) being a complex of constant sheaves implies that the cohomology $H^k\Big(\sL_{\underline A}\boxtimes \bar j_{\underline F,\,\underline A}^{-1}\sW_i(\sL_{\underline A})\Big)=\sL_{\underline A}\otimes_\C H^k (\bar j_{\underline F,\,\underline A}^{-1}\sW_i(\sL_{\underline A}))$ of the right hand side is a constant sheaf for each $i$. The cohomology $H^k (\bar j_{\underline F,\,\underline A}^{-1}\sW_i(\sL_{\underline A}))$ is a $\C$-vector space, so $\sL_{\underline A}$ being not a constant sheaf implies that $H^k (\bar j_{\underline F,\,\underline A}^{-1}\sW_i(\sL_{\underline A}))$ vanishes for each $k$. Hence $\bar j_{\underline F,\,\underline A}^{-1}\sW_i(\sL_{\underline A})=0$, and then by (\ref{qiao}) we have
$\bar j_{F,\,A}^{-1}\sW_i(\sL_A)=0.$ Using the same method, $\bar j_{F,\,A}^{-1}(\bar k_i^A)_{!*}(\bar \ell_i^A)_!(\sL_A^{-1}[d_A])=0$, where $\sL_A^{-1}={\rm Hom}_{\C_{T_A^{\rm an}}}(\sL_A,\C_{T_A^{\rm an}})$. Applying the Verdier duality functor, we have
$
\bar j_{F,\,A}^!\sW_i(\sL_A)=\bar j_{F,\,A}^!(\bar k_i^A)_{!*}(\bar \ell_i^A)_*(\sL_A[d_A])=0.
$
\end{proof}

\begin{rem}
It should be noted that a $D$-module version of Lemma \ref{1214} was given by A. Steiner in \cite[Lemma 9.1 (b), (d)]{AS}. His result and mine are equivalent via the Riemann-Hilbert correspondence.
\end{rem}

For the construction of the canonical epimorphism in Theorem \ref{yaode}, we need the following lemma.

\begin{lem}\label{1215}
For any $F\in \mathfrak F_i(A)$, let $h_F\colon T_F\to U_i(A)$ be the natural closed immersion. Then
\begin{eqnarray*}
(\bar k_i^A)^{-1}(\sW_i(\sL_A)/\sW_{i-1}(\sL_A))=\bigoplus_{\substack{F\in \mathfrak F_i(A)\\\sL_F\in\pi_{F,A}^{\diamond}(\sL_A)}}(h_F)_*(\sL_F[d_F]).
\end{eqnarray*}
\end{lem}

\begin{proof}
Denote by $j_i^A$ the open immersion $U_{i-1}(A)\to U_i(A)$. By the proof of \cite[Proposition 8.2.11]{HTT}, we have
\begin{eqnarray*}\label{422}
(\bar k_i^A)^{-1}\sW_{i-1}(\sL_A)=(j_i^A)_{!*}(\bar \ell_{i-1}^A)_*(\sL_A[d_A])=\tau_{<i-d_A}(\bar \ell_i^A)_*(\sL_A[d_A]).
\end{eqnarray*}
Here, for any complex $\sF^\bullet$ of sheaves on a scheme and any $i\in\Z$, we use the notation that
\begin{align*}
&\tau^{<i}\sF^\bullet=\Big(\cdots\sF^{i-4}\to\sF^{i-3}\to\sF^{i-2}\to{\rm ker}(\sF^{i-1}\to\sF^i)\to0\to0\to\cdots\Big);\\
&\tau^{\geq i}\sF^\bullet=\Big(\cdots\to0\to0\to{\rm coker}(\sF^{i-1}\to\sF^i)\to\sF^{i+1}\to\sF^{i+2}\to\cdots\Big).\end{align*}
Since $\bar k_i^A$ is an open immersion, then
$$(\bar k_i^A)^{-1}(\sW_i(\sL_A)/\sW_{i-1}(\sL_A))=(\bar k_i^A)^{-1}(\sW_i(\sL_A))/(\bar k_i^A)^{-1}(\sW_{i-1}(\sL_A))=\tau^{\geq i-d_A}(\bar \ell_i^A)_*(\sL_A[d_A])$$ is a perverse sheaf on $U_i(A)$ supported on $U_i(A)-U_{i-1}(A)=\mathop\bigsqcup\limits_{F\in \mathfrak F_i(A)} h_F(T_F)$.
As a result,
\begin{eqnarray}\label{18326}(\bar k_i^A)^{-1}(\sW_i(\sL_A)/\sW_{i-1}(\sL_A))=\bigoplus_{F\in \mathfrak F_i(A)}\tau^{\geq i-d_A}(h_F)_*h_F^{-1}(\bar \ell_i^A)_*(\sL_A[d_A]).
\end{eqnarray}
First assume that $A$ is normal. By (\ref{ji4}) and (\ref{18326}), we therefore have
\begin{align*}
&\quad\;(\bar k_i^A)^{-1}(\sW_i(\sL_A)/\sW_{i-1}(\sL_A))\\
&= \bigoplus_{F\in \mathfrak F_i(A)}\tau^{\geq i-d_A}(h_F)_*h_F^{-1}(\bar \ell_i^A)_*(\sL_A[d_A])\\
&=\bigoplus_{F\in \mathfrak F_i(A)}(h_F)_*\tau^{\geq i-d_A}h_F^{-1}(\bar k_i^A)^{-1}(\bar k_i^A)_*(\bar \ell_i^A)_*(\sL_A[d_A])\\
&= \bigoplus_{F\in \mathfrak F_i(A)}(h_F)_*\tau^{\geq i-d_A}\bar j_{F,\,A}^{-1}(\bar j_A)_*(\sL_A[d_A])\\
&= \bigoplus_{\substack{F\in \mathfrak F_i(A)\\\sL_F\in\pi_{F,A}^{\diamond}(\sL_A)}}(h_F)_*\tau^{\geq i-d_A}\Big(\sL_F[d_F]\otimes_\C \textstyle\bigwedge^\bullet\C^{d_A-d_F}[d_A-d_F]\Big)\\
&= \bigoplus_{\substack{F\in \mathfrak F_i(A)\\\sL_F\in\pi_{F,A}^{\diamond}(\sL_A)}}(h_F)_*(\sL_F[d_F]).
\end{align*}
This proves the lemma for any normal matrix $A$.

For general $A$, choose $a_{N+1},\ldots,a_{N+\ell}\in\Z A$ such that $\mathbb R_{\geq0}A\cap\Z A=\sum\limits_{j=1}^{N+\ell}\N a_j$. Thus $\widetilde A:=(a_1,\ldots,a_{N+\ell})$ is a normal matrix and $\sL_A$ can be viewed as a local system $\sL_{\widetilde A}$ on $T_{\widetilde A}=T_A$. The inclusion $A\subset\widetilde A$ induces two finite morphisms $\pi:X_{\widetilde A}\to X_A$ and $\pi_i:U_i(\widetilde A)\to U_i(A)$. Then
\begin{align*}\begin{aligned}
&\quad(\bar k_i^A)^{-1}(\sW_i(\sL_A)/\sW_{i-1}(\sL_A))\\
&=(\bar k_i^A)^{-1}\pi_*(\sW_i(\sL_{\widetilde A})/\sW_{i-1}(\sL_{\widetilde A}))\\
&=(\pi_i)_*(\bar k_i^{\widetilde A})^{-1}(\sW_i(\sL_{\widetilde A})/\sW_{i-1}(\sL_{\widetilde A}))\\
&=\bigoplus_{\substack{\widetilde F\in\mathfrak F_i(\widetilde A)\\\sL_{\widetilde F}\in\pi_{\widetilde F,\widetilde A}^{\diamond}(\sL_{\widetilde A})}}(\pi_i)_*(h_{\widetilde F})_*(\sL_{\widetilde F}[d_{\widetilde F}])\\
&=\bigoplus_{\substack{\widetilde F\in\mathfrak F_i(\widetilde A)\\\sL_{\widetilde F}\in\pi_{\widetilde F,\widetilde A}^{\diamond}(\sL_{\widetilde A})}}(h_F)_*(\pi_{F,\,\widetilde F})_*(\sL_{\widetilde F}[d_{\widetilde F}])\\
&=\bigoplus_{\substack{\widetilde F\in\mathfrak F_i(\widetilde A)\\\sL_{\widetilde F}\in\pi_{\widetilde F,\widetilde A}^{\diamond}(\sL_{\widetilde A})}}\Big(\bigoplus_{\sL_F\in\pi_{F,\widetilde F}^{\diamond}(\sL_{\widetilde F})}(h_F)_*(\sL_F[d_F])\Big)\\
&=\bigoplus_{\substack{F\in \mathfrak F_i(A)\\\sL_F\in\pi_{F,A}^{\diamond}(\sL_A)}}(h_F)_*(\sL_F[d_F]),
\end{aligned}\end{align*}
where the first and fifth equalities use Lemma \ref{4134}, the second uses the proper base change theorem, the third is the case that we have just proved, the fourth holds by $\pi_i\circ h_{\widetilde F}=h_F\circ\pi_{F,\widetilde F}:T_{\widetilde F}\to U_i(A)$ and the last follows from the fact that $\pi_{F,A}^{\diamond}(\sL_A)=\pi_{F,\widetilde F}^{\diamond}(\pi_{\widetilde F,\widetilde A}^{\diamond}(\sL_{\widetilde A}))$.
\end{proof}

\subsection{Proof of Theorem \ref{yaode}.}
\begin{proof}[Proof of Theorem \ref{yaode}]
By \cite[Corollary 8.29]{HTT}, the minimal extension functor $(\bar k_i^A)_{!*}\colon{\rm Perv}(U_i(A))\to{\rm Perv}(X_A)$ preserves injectivity and surjectivity. Applying $(\bar k_i^A)_{!*}$ to the short exact sequence
$$0\to(j_i^A)_{!*}(\bar \ell_{i-1}^A)_*(\sL_A[d_A])\to(\bar \ell_i^A)_*(\sL_A[d_A])\to\bigoplus_{\substack{F\in \mathfrak F_i(A)\\\sL_F\in\pi_{F,A}^{\diamond}(\sL_A)}}(h_F)_*(\sL_F[d_F])\to0$$
of perverse sheaves on $U_i(A)$ given by Lemma \ref{1215}, then the composition
$$\sW_{i-1}(\sL_A)\hookrightarrow\sW_i(\sL_A)\twoheadrightarrow\bigoplus_{\substack{F\in \mathfrak F_i(A)\\\sL_F\in\pi_{F,A}^{\diamond}(\sL_A)}}(\bar j_{F,\,A})_{!*}(\sL_F[d_F])$$
is trivial, and hence we have a canonical epimorphism
\begin{eqnarray*}\label{mei}
\alpha_{\sL_A}\colon \sW_i(\sL_A)/\sW_{i-1}(\sL_A)\twoheadrightarrow\bigoplus_{\substack{F\in \mathfrak F_i(A)\\\sL_F\in\pi_{F,A}^{\diamond}(\sL_A)}}(\bar j_{F,\,A})_{!*}(\sL_F[d_F]).
\end{eqnarray*}
This proves part (1) of Theorem \ref{yaode}.

It remains to prove part (2) of Theorem \ref{yaode}. We do this in 5 steps. In Step (i), we prove part (2) of Theorem 1.7 for $A=(1)$. In Step (ii), we prove the theorem for the $n\times n$ identity matrix $I_n$. In Step (iii), we prove the theorem when $\mathbb R_{\geq0}A$ is a simplicial cone. In Step (iv), we prove the theorem for a normal matrix $A$. In Step (v), we prove the theorem for general $A$.

\underline{Step (i).} Part (2) of Theorem \ref{yaode} holds if $A=I_1=(1)$.

In this case, $\bar j_A$ is the open immersion $\G_{\rm m}\hookrightarrow\A^1$, $U_0(I_1)=\G_{\rm m}$ and $U_1(I_1)=\A^1$.
If $\sL_A$ is not the constant sheaf, then $\sW_0(\sL_A)=\sW_1(\sL_A)=(\bar j_A)_{!*}(\sL_A[1])$. Otherwise, $\sW_1(\sL_A)/\sW_0(\sL_A)=(\bar j_A)_*(\C_{\G_{\rm m}^{\rm an}}[1])/(\bar j_A)_{!*}(\C_{\G_{\rm m}^{\rm an}}[1])=[0]_*\C$, where $[0]$ is the complement of $\bar j_A$. This proves Step (i).

\underline{Step (ii).} Part (2) of Theorem \ref{yaode} holds if $A$ is the $n\times n$ identity matrix $I_n$.

In this case, $\bar j_A$ is the open immersion $\G_{\rm m}^n\hookrightarrow \A^n$ and $\sL_A=\sL_1\boxtimes\cdots\boxtimes\sL_n$ for some local systems $\sL_1,\ldots,\sL_n$ on $\G_{\rm m}$. For any $0\leq i_1,\ldots,i_n\leq 1$ with $\sum\limits_{t=1}^ni_t=i$, $\bar j_A$ factors as
$$\G_m^n\x{\bar \ell_{i_1}^{I_1}\times\cdots\times\bar \ell_{i_n}^{I_1}}U_{i_1}(I_1)\times\cdots\times U_{i_n}(I_1)\x{k}U_i(I_n)\x{\bar k_i^A}\A^n.$$
The injectivity of the functor $(\bar k_i^A)_{!*}$ on perverse sheaves shows that
\begin{eqnarray*}
&&\mathop{\boxtimes}\limits_{t=1}^n\sW_{i_t}(\sL_t)= \mathop{\boxtimes}\limits_{t=1}^n(\bar k_{i_t}^{I_1})_{!*}(\bar\ell_{i_t}^{I_1})_*(\sL_{t}[1])=(\bar k_i^A)_{!*}k_{!*}(\prod\limits\limits_{t=1}^n\bar\ell_{i_t}^{I_1})_*(\sL_A[n])\\
&\subseteq&(\bar k_i^A)_{!*}k_{*}(\prod\limits\limits_{t=1}^n\bar\ell_{i_t}^{I_1})_*(\sL_A[n])=(\bar k_i^A)_{!*}(\bar\ell_i^{A})_*(\sL_A[n])=\sW_i(\sL_A).
\end{eqnarray*}
Since the restriction of $\mathop{\boxtimes}\limits_{t=1}^n\sW_{i_t}(\sL_t)$ and $\sW_i(\sL_A)$ on $U_{i_1}(I_1)\times\cdots\times U_{i_n}(I_1)$ all coincide with $\mathop{\boxtimes}\limits_{t=1}^n(\bar\ell_{i_t}^{I_1})_*(\sL_{t}[1])$, then
$\sum\limits\limits_{\substack{0\leq i_t\leq 1\\i_1+\cdots+i_n=i}}(\mathop{\boxtimes}\limits_{t=1}^n\sW_{i_t}(\sL_t))$ is a perverse subsheaf of $\sW_i(\sL_A)$, whose restriction on
$\mathop{\bigcup}\limits_{\substack{0\leq i_t\leq 1\\i_1+\cdots+i_n=i}}(U_{i_1}(I_1)\times\cdots\times U_{i_n}(I_1))=U_i(I_n)$ is $(\bar\ell_i^A)_*(\sL_A[n])$. By the definition of the minimal extension of perverse sheaves, $\sW_i(\sL_A)=(\bar k_i^A)_{!*}(\bar\ell_i^A)_*(\sL_A[n])$ is the smallest perverse subsheaf of $\sW_n(\sL_A)$ whose restriction on $U_i(I_n)$ coincides with $(\bar\ell_i^A)_*(\sL_A[n])$. Consequently,
$$\sW_i(\sL_A)=\sum\limits_{\substack{0\leq i_t\leq 1\\i_1+\cdots+i_n=i}}\sW_{i_1}(\sL_1)\boxtimes\cdots\boxtimes\sW_{i_n}(\sL_n).$$
Set $\sW_{-1}(\sL_i)=0$. By the exactness of the bi-functor $\boxtimes$ on perverse sheaves, we obtain
$$\sW_i(\sL_A)/\sW_{i-1}(\sL_A)=\bigoplus_{\substack{0\leq i_t\leq 1\\i_1+\cdots+i_n=i}}(\sW_{i_1}(\sL_1)/\sW_{i_1-1}(\sL_1))\boxtimes\cdots\boxtimes(\sW_{i_n}(\sL_n)/\sW_{i_n-1}(\sL_n)).$$
This proves part (2) of Theorem \ref{yaode} for $A=I_n$, and hence for any matrix $A$ with $d_A=N$.

\underline{Step (iii).} Part (2) of Theorem \ref{yaode} holds if $\mathbb R_{\geq0}A$ is a simplicial cone.

In this case, there is a subset $A_0$ of $A$ with $d_A$ elements such that $\mathbb R_{\geq0}A=\mathbb R_{\geq0}A_0$. Then $A_0\subset A$ induces a finite and surjective morphism $p\colon X_A\to X_{A_0}$ by Lemma \ref{4132}. The map $F\mapsto F_0:= F\cap A_0$ defines a bijection between the set of faces of $A$ and that of $A_0$. We have
\begin{align*}
&\quad p_*(\sW_i(\sL_A)/\sW_{i-1}(\sL_A))\\
&=\bigoplus_{\sL_{A_0}\in\pi_{A_0,A}^{\diamond}(\sL_A)}\sW_i(\sL_{A_0})/\sW_{i-1}(\sL_{A_0})\\
&=\bigoplus_{\sL_{A_0}\in\pi_{A_0,A}^{\diamond}(\sL_A)}\Big(\bigoplus_{\substack{F_0\in \mathfrak F_i(A_0)\\\sL_{F_0}\in\pi_{F_0, A_0}^{\diamond}(\sL_{A_0})}}(\bar j_{F_0,\,A_0})_{!*}(\sL_{F_0}[d_{F_0}])\Big)\\
&=\bigoplus_{\substack{F\in \mathfrak F_i(A)\\\sL_F\in\pi_{F, A}^{\diamond}(\sL_A)}}\Big(\bigoplus_{\sL_{F_0}\in\pi_{F_0,F}^{\diamond}(\sL_F)}(\bar j_{F_0,\,A_0})_{!*}(\sL_{F_0}[d_{F_0}])\Big)\\
&=\bigoplus_{\substack{F\in \mathfrak F_i(A)\\\sL_F\in\pi_{F, A}^{\diamond}(\sL_A)}}(\bar j_{F_0,\,A_0})_{!*}(\pi_{F_0,\,F})_*(\sL_{F}[d_{F}])\\
&=\bigoplus_{\substack{F\in \mathfrak F_i(A)\\\sL_F\in\pi_{F, A}^{\diamond}(\sL_A)}}\pi_*(\bar j_{F,\,A})_{!*}(\sL_{F}[d_{F}]).
\end{align*}
Let's give a brief explaination of the above equation. The first and the fourth equalities use Lemma \ref{4134}, the second applies Step (ii) to $A_0$, the third follows from the fact that $\pi_{F_0,A}^{\diamond}(\sL_A)=\pi_{F_0,F}^{\diamond}(\pi_{F,A}^{\diamond}(\sL_A))=\pi_{F_0,A_0}^{\diamond}(\pi_{A_0,A}^{\diamond}(\sL_A))$ and the last uses part (2) of Lemma \ref{4131}.
This proves that $\pi_*(\alpha_{\sL_A})$ is an isomorphism. Hence, so is $\alpha_{\sL_A}$ by the surjectivity of $\pi$. This proves Step (iii).

\underline{Step (iv).} Part (2) of Theorem \ref{yaode} holds if $A$ is normal.

For this step, assume that $A$ is normal. Recall that $F_1,\ldots,F_k$ are all facet of $A$ such that $\sL_A$ is the inverse image of some local system on $T_{F_i}$. Choose a section $\iota$ of the epimorphism $\tau\colon \Z A\to \Z A/\Z F_0$ where $F_0=F_1\cap\cdots\cap F_k$, and set $\epsilon_{F_0}=\sum\limits_{a_j\in F_0}a_j$. By the map
$F\mapsto\overline F:=(F,-\epsilon_{F_0})\mapsto\underline F:=\iota(\tau(F)),$
the three sets \{faces of $A$ containing $F_0$\}, \{faces of $\overline A$\} and \{faces of $\underline A$\} are in bijective correspondence. As $T_A = T_{\overline A}$, the local system $\sL_A$ on $T_A$ can be viewed as a local system $\sL_{\overline A}$ on $T_{\overline A}$. There exists a local system $\sL_{F_0}$ on $T_{F_0}$ and a local system $\sL_{\underline A}$ on $T_{\underline A}$ such that $\sL_{\overline A} = \sL_{F_0}\boxtimes\sL_{\underline A}$ on $T_{\overline A}= T_{F_0}\times T_{\underline A }$ defined by $\iota$. The assumption on $F_1,\ldots,F_k$ implies that $\mathbb R_{\geq0}\underline A$ is a simplicial cone. Then
\begin{align}\begin{aligned}\label{183179}
&\quad\;\sW_i(\sL_{\overline A})/\sW_{i-1}(\sL_{\overline A})\\
&=\sL_{F_0}[d_{F_0}]\boxtimes(\sW_i(\sL_{\underline A})/\sW_{i-1}(\sL_{\underline A}))\\
&=\sL_{F_0}[d_{F_0}]\boxtimes\bigoplus_{\substack{\underline F\in \mathfrak F_i(\underline A)\\\sL_{\underline F}\in\pi^{\diamond}_{\underline F,\underline A}(\sL_{\underline A})}}(\bar j_{\underline F,\,\underline A})_{!*}(\sL_{\underline F}[d_{\underline F}])\\
&=\bigoplus_{\substack{\underline F\in \mathfrak F_i(\underline A)\\\sL_{\underline F}\in\pi^{\diamond}_{\underline F,\underline A}(\sL_{\underline A})}}({\rm id}_{T_{F_0}}\times \bar j_{\underline F,\,\underline A})_{!*}(\sL_{F_0}[d_{F_0}]\boxtimes\sL_{\underline F}[d_{\underline F}])\\
&=\bigoplus_{\substack{\overline F\in \mathfrak F_i(\overline A)\\\sL_{\overline F}\in\pi_{\overline F,\overline A}^{\diamond}(\sL_{\overline A})}}(\bar j_{\overline F,\,\overline A})_{!*}(\sL_{\overline F}[d_{\overline F}]),
\end{aligned}\end{align}
where the first equality uses (\ref{wen}), the second uses Step (iii), the third follows from the last paragraph of $\S 1.2$ and the last is trivial.

For any face $F$ of $A$ not containing $F_0$, there exists a facet $F'$ of $A$ different from
$F_1,\ldots,F_k$ such that $F\prec F'$. So $\sL_A$ is not the inverse image under $\pi_{F,A}$ of any local system on $T_F$. By (5.2), $\bar j^{-1}_{F,A}\sW_i(\sL_A) =\bar j^!_{F,A}\sW_i(\sL_A) = 0$, and hence $\zeta^{-1}\sW_i(\sL_A) = \zeta^!\sW_i(\sL_A) = 0$, where $\zeta\colon\mathop{\bigcup}\limits_{F_0\nprec F\prec A} T_F\to X_A$ is the complement of the open immersion $j \colon  X_{\overline A}\to X_A$.
Then the canonical isomorphism $j^{-1}\sW_i(\sL_A)\simeq\sW_i(\sL_{\overline A})$ induces two canonical isomorphisms
\begin{eqnarray}\label{12151}j_!\sW_i(\sL_{\overline A})\simeq\sW_i(\sL_A)\simeq j_*\sW_i(\sL_{\overline A}).\end{eqnarray}
By the proof of Lemma \ref{4131}, there are two exact functors $j_*,\,j_!:{\rm Perv}(X_{\overline A})\to{\rm Perv}(X_A)$. We thus have a commutative diagram
\begin{eqnarray}\begin{split}\label{12152}\xymatrix{0\ar[r]&j_!\sW_{i-1}(\sL_{\overline A})\ar[r]\ar[d]^\simeq&j_!\sW_{i}(\sL_{\overline A})\ar[r]\ar[d]^\simeq&j_!(\sW_{i}(\sL_{\overline A})/\sW_{i-1}(\sL_{\overline A}))\ar[r]\ar[d]^\simeq&0\\
0\ar[r]&j_*\sW_{i-1}(\sL_{\overline A})\ar[r]&j_*\sW_{i}(\sL_{\overline A})\ar[r]&j_*(\sW_{i}(\sL_{\overline A})/\sW_{i-1}(\sL_{\overline A}))\ar[r]&0}
\end{split}\end{eqnarray}
of perverse sheaves with exact rows. So we obtain
\begin{eqnarray*}
\sW_i(\sL_A)/\sW_{i-1}(\sL_A)&=&j_{!*}(\sW_i(\sL_{\overline A})/\sW_{i-1}(\sL_{\overline A}))\\
&=&\bigoplus_{\substack{\overline F\in \mathfrak F_i(\overline A)\\\sL_{\overline F}\in\pi_{\overline F,\overline A}^{\diamond}(\sL_{\overline A})}}j_{!*}(\bar j_{\overline {F},\,\overline A})_{!*}(\sL_{\overline {F}}[d_{\overline {F}}])\\
&=&\bigoplus_{\substack{F\in \mathfrak F_i(A)\\\sL_F\in\pi_{F,A}^{\diamond}(\sL_A)}}(\bar j_{F,\,A})_{!*}(\sL_F[d_F]),
\end{eqnarray*}
where the first equality uses (\ref{12151}) and (\ref{12152}), the second uses (\ref{183179}) and the last is obvious.
This proves Theorem part (2) of \ref{yaode} when $A$ is normal.

\underline{Step (v).} Part (2) of Theorem \ref{yaode} holds for general $A$.

 For general $A$, choose $a_{N+1},\ldots,a_{N+\ell}\in\Z A$ such that $\mathbb R_{\geq 0}A\cap\Z A=\sum\limits_{j=1}^{N+\ell}\N a_j$, and let $\widetilde A=(A,a_{N+1},\ldots,a_{N+\ell})$. The map $\widetilde F\mapsto F:=\widetilde F\cap A$ defines a bijection between the set of faces of $\widetilde A$ and that of $A$. The local system $\sL_A$ on $T_A=T_{\widetilde A}$ can be viewed as a local system $\sL_{\widetilde A}$ on $T_{\widetilde A}$. Let $\pi\colon X_{\widetilde A}\to X_A$ be the morphism induced by $A\subset \widetilde A$.
Using the results of Step (iv) and repeating the computation in Step (iii), we also have
\begin{align*}
&\quad\;\sW_i(\sL_A)/\sW_{i-1}(\sL_A)\\
&=\pi_*(\sW_i(\sL_{\widetilde A})/\sW_{i-1}(\sL_{\widetilde A}))\\
&=\bigoplus_{\substack{\widetilde F\in \mathfrak F_i(\widetilde A)\\\sL_{\widetilde F}\in\pi_{\widetilde F,\widetilde A}^{\diamond}(\sL_{\widetilde A})}}\pi_*(\bar j_{\widetilde F,\,\widetilde A})_{!*}(\sL_{\widetilde F}[d_{\widetilde F}])\\
&=\bigoplus_{\substack{\widetilde F\in \mathfrak F_i(\widetilde A)\\\sL_{\widetilde F}\in\pi_{\widetilde F,\widetilde A}^{\diamond}(\sL_{\widetilde A})}}(\bar j_{F,\,A})_{!*}(\pi_{F,\,\widetilde F})_*(\sL_{\widetilde F}[d_{\widetilde F}])\\
&=\bigoplus_{\substack{\widetilde F\in \mathfrak F_i(\widetilde A)\\\sL_{\widetilde F}\in\pi_{\widetilde F,\widetilde A}^{\diamond}(\sL_{\widetilde A})}}\Big(\bigoplus_{\sL_F\in\pi_{F,\widetilde F}^{\diamond}(\sL_{\widetilde F})}(\bar j_{F,\,A})_{!*}(\sL_{F}[d_{F}])\Big)\\
&=\bigoplus_{\substack{ F\in \mathfrak F_i( A)\\\sL_{ F}\in\pi_{ F, A}^{\diamond}(\sL_{ A})}}(\bar j_{F,\,A})_{!*}(\sL_{F}[d_{F}]).
\end{align*}
This completes the proof of Theorem \ref{yaode}.
\end{proof}

\begin{rem}
Theorem \ref{yaode} is a perverse sheaf analog of Theorem 1.4, but it cannot be proven from Theorem 1.4.

For example, let $A=\left(\begin{array}[c]{lll}1&1&1\\0&1&2\end{array}\right)$, $\beta=\left(\begin{array}[c]{l}\frac12\\0\end{array}\right)$, and let $\sL_A$ be the local system on $T_A={\rm Spec}\,\C[t_1^{\pm1},t_2^{\pm1}]$ defined by the multi-valued function $t_1^{\frac12}$. We have $\sL_A[2]={\rm DR}_{T_A}(\sO_{T_A}^{\beta'})$ if and only if $\beta'\in\beta+\Z^2$. If we can reduce Theorem 1.7 to Theorem 1.4 by the Riemann-Hilbert correspondence,
there must exists $\beta'\in\beta+\Z^2$ such that $W_i(A,\beta')=\widetilde W_i(A,\beta')$ for any $i$. By Theorem \ref{baaa}, the holonomic $D_A$-module $W_1(A,\beta')/W_0(A,\beta')=\mathop{\bigoplus}\limits_{\substack{F\in\mathfrak F_1(A)\\\beta'\in\C F}}(j_{F,A})_{!+}\sO_{T_F}^{\beta'}$ has length $\leq1$. By Theorem \ref{yaode}, $\widetilde W_1(A,\beta')/\widetilde W_{0}(A,\beta')=(j_{F_1,A})_{!+}\sO_{T_{F_1}}^{\beta}\oplus(j_{F_2,A})_{!+}\sO_{T_{F_2}}^{\beta''}$ has length 2, where $\beta''=(\frac{1}{2},1)^{\rm t}$. These contradict to the fact that $W_i(A,\beta')=\widetilde W_i(A,\beta')$ for any $i$.

\end{rem}
\subsection{Counterexample to Theorem \ref{yaode} in the non-simplicial case}
The following example shows that $\alpha_{\sL_A}$ may not be an isomorphism if $A$ is not simplicial relative to $\sL_A$.

\if<\begin{exam}\label{hui}
Let \[A=\left(\begin{array}[c]{llll}1&0&0&1\\0&1&0&1\\0&0&1&-1\end{array}\right).\]
For any $1\leq i\leq 3$, denote by $\sK_i$ the kernel of the canonical epimorphism
$$\alpha_i\colon \sW_i(\C_{T_A^{\rm an}})/\sW_{i-1}(\C_{T_A^{\rm an}})\twoheadrightarrow\bigoplus_{G\in \mathfrak F_i(A)}(\bar j_{G,A})_{!*}(\C_{T_{G}^{\rm an}}[d_G]).$$
Then $\sK_1=\sK_3=0$ and $\sK_2=(i_0)_*\C$, where $i_0\colon {\rm Spec}\,\C\to X_A$ be the closed immersion defined by the unique $T_A$-fixed point $x_0$ on $X_A(\C)$.
\end{exam}

\begin{proof}
For each $i$, put $\sQ_i=\sW_i(\C_{T_A^{\rm an}})/\sW_{i-1}(\C_{T_A^{\rm an}})$. By \cite[8.2.11]{HTT},
\begin{align}\label{s3}&\sW_0(\C_{T_A^{\rm an}})=\tau_{<0}(\bar k_2^A)_*\tau_{<-1}(j_1)_*\tau_{<-2}(\bar \ell_1^A)_*(\C_{T_A^{\rm an}}[3]);\\
\label{s1}&\sW_1(\C_{T_A^{\rm an}})=\tau_{<0}(\bar k_2^A)_*\tau_{<-1}(\bar \ell_2^A)_*(\C_{T_A^{\rm an}}[3]);\\
\label{s2}&\sW_2(\C_{T_A^{\rm an}})=\tau_{<0}(\bar j_A)_*(\C_{T_A^{\rm an}}[3]),\end{align}
where $j_i$ is the open immersion $U_1(A)\hookrightarrow U_2(A)$. This proves that $H^i\sQ=0$ for any $i\geq0$.
By Lemma \ref{1215}, $(\bar k_i^A)^{-1}(\alpha_i)$ is an isomorphism, and the kernel $\sK_i$ of $\alpha_i$ is therefore a perverse sheaf on $X_A$ supported on $X_A-U_i(A)$. In particular, $\sK_3=0$, $\sK_2$ is supported on $\{x_0\}$ and $\sK_1$ is supported on $X_A-U_1(A)$.

For any $F\in\mathfrak F_2(A)$, we have a distinguished triangle
\begin{align}\label{w3}
\bar j_{F,A}^{-1}\sK_1\to\bar j_{F,A}^{-1}\sQ_1\to\bigoplus_{G\in \mathfrak F_1(A)}\bar j_{F,A}^{-1}(\bar j_{G,A})_{!*}(\C_{T_{G}^{\rm an}}[2]).
\end{align}
Since $\sK_1$ is a perverse on $X_A$ supported on $X_A-U_1(A)=\mathop{\bigcup}\limits_{\substack{F\prec A,\,d_F\leq1}}T_F$, and $T_F\hookrightarrow X_A-U_1(A)$ is an open immersion, then $\bar j_{F,A}^{-1}\sK_1$ is a perverse sheaves on $T_F$. For any facet $G$ of $A$, $\bar j_G:T_G\to X_G$ can be identified with the inclusion $\G^2_{\rm m}\hookrightarrow\A^2$. So we have
\begin{align}\begin{aligned}\label{w2}
&\quad\;\bigoplus_{G\in \mathfrak F_1(A)}\bar j_{F,A}^{-1}(\bar j_{G,A})_{!*}(\C_{T_{G}^{\rm an}}[2])\\
&=\bigoplus_{G\in \mathfrak F_1(A)}\bar j_{F,A}^{-1}(\bar i_{G,A})_*(\bar j_G)_{!*}(\C_{T_{G}^{\rm an}}[2])\\
&=\bigoplus_{G\in \mathfrak F_1(A)}\bar j_{F,A}^{-1}(\bar i_{G,A})_*(\C_{X_{G}^{\rm an}}[2])\\
&=\C^2_{T_{F}^{\rm an}}[2].
\end{aligned}\end{align}
Applying Theorem \ref{yaode} to the $2\times 2$-matrix $I_2$ and the constant sheaf $\C_{T_{I_2}^{\rm an}}$ on $T_{I_2}=\G^2_{\rm m}$, we have
\begin{align*}
\iota_0^{-1}(\sW_1(\C_{T_{I_2}^{\rm an}})/\sW_0(\C_{T_{I_2}^{\rm an}}))=\bigoplus_{G\in \mathfrak F_1(I_2)}\iota_0^{-1}(\bar j_{G,I_2})_{!*}(\C_{T_{G}^{\rm an}}[1])=\C^2[1],
\end{align*}
where $\iota_0:{\rm Spec}\,\C\to\A^2$ is the morphism defined by $0\in\A^2(\C)$.
By the proof (\ref{qiao}) in Lemma \ref{1214}, $j_{\underline F,\,\underline A}$ for $F\in\mathfrak F_2(A)$ can be identified with the morphism $\iota_0$. So
\begin{align*}
 \bar j_{F,A}^{-1}\sW_i(\C_{T_A^{\rm an}})=\C_{T_F^{\rm an}}[1]\otimes\iota_0^{-1}\sW_i(\C_{T_{I_2}^{\rm an}}),
\end{align*}
and hence
\begin{align}\label{w1}
 \bar j_{F,A}^{-1}(\sW_1(\C_{T_A^{\rm an}})/\sW_{0}(\C_{T_A^{\rm an}}))=\C_{T_F^{\rm an}}[1]\otimes\iota_0^{-1}(\sW_1(\C_{T_{I_2}^{\rm an}})/\sW_0(\C_{T_{I_2}^{\rm an}}))=\C^2_{T_F^{\rm an}}[2].
\end{align}
Substituting (\ref{w2}) and (\ref{w1}) to the distinguished triangle (\ref{w3}), we have a distinguished triangle
\begin{align}\label{w4}
  \bar j_{F,A}^{-1}\sK_1\to\C^2_{T_{F}^{\rm an}}[2]\to \C^2_{T_{F}^{\rm an}}[2],
\end{align}
and an exact sequence
$$H^{-2}\bar j_{F,A}^{-1}\sK_1\to\C^2_{T_{F}^{\rm an}}\to \C^2_{T_{F}^{\rm an}}.$$
Since $\bar j_{F,A}^{-1}\sK_1$ is a perverse sheaf on $T_F$, then $H^{-2}\bar j_{F,A}^{-1}\sK_1=0$ and hence the second morphism in (\ref{w4}) is a quasi-isomorphism. This proves that $\bar j_{F,A}^{-1}\sK_1=0$ for any $F\in\mathfrak F_2(A)$. So $\sK_1$ is supported on $\{x_0\}$.

For any $k$, $\bar j_{F_k}$ is the open immersion $T_{F_k}={\rm Spec}\,\C[x_k^\pm]\hookrightarrow X_{F_k}={\rm Spec}\,\C[x_k]$. Then
\begin{align}\begin{aligned}
&\quad\;(\bar j_{F_k,\,A})_{!*}(\C_{T_{F_k}^{\rm an}}[1])\\
&=(\bar i_{F_k,\,A})_*(\bar j_{F_k})_{!*}(\C_{T_{F_k}^{\rm an}}[1])=(\bar i_{F_k,\,A})_*\C_{X_{F_k}^{\rm an}}[1].
\end{aligned}\end{align}

Consequently, $i_0^{-1}(\bar j_{F_k,\,A})_{!*}(\C_{T_{F_k}^{\rm an}}[1])=\C[1]. $ Thus, according to the short exact sequence $$0\to\sK\to\sQ\to\bigoplus_{k=1}^4(\bar j_{F_k,A})_{!*}(\C_{T_{F_k}^{\rm an}}[1])\to0$$ of perverse sheaves, we have $H^i\sQ=0$ for any $i\neq-1$.
Applying $H^i$ to the distinguished triangle $$i_0^{-1}\sK\to i_0^{-1}\sQ\to\bigoplus_{k=1}^4i_0^{-1}(\bar i_{F_k,\,A})_{*}\C_{X_{F_k}}[1],$$ we thus have a short exact sequence of $\C$-vector spaces:
\begin{eqnarray}\label{1911}0\to H^{-1}i_0^{-1}\sQ\to\C^4\to H^0i_0^{-1}\sK\to0.\end{eqnarray}
By (\ref{diao}) and (\ref{s2}), $$H^{-1}i_0^{-1}\sW_2(\C_{T_A^{\rm an}})= H^{-1}i_0^{-1}(\bar j_A)_*(\C_{T_A^{\rm an}}[3])=H^2i_0^{-1}(\bar j_A)_*(\C_{T_A^{\rm an}})=H^2(\textstyle\bigwedge^\bullet\C^3)=\C^3.$$
Applying $H^i$ to the distinguished triangle
$$i_0^{-1}\tau_{<0}(\bar k_2^A)_*\tau_{<-1}(\bar \ell_2^A)_*(\C_{T_A^{\rm an}}[3])\to i_0^{-1}\tau_{<0}(\bar j_A)_*(\C_{T_A^{\rm an}}[3])\to i_0^{-1}\sQ,$$
we obtain an exact sequence
\begin{eqnarray}\label{1912}0\to H^{-1}i_0^{-1}\sW_1(\C_{T_A^{\rm an}})\to H^{-1}i_0^{-1}\sW_2(\C_{T_A^{\rm an}})=\C^3\to H^{-1}i_0^{-1}\sQ\to0.\end{eqnarray}
Combining (\ref{1911}) with (\ref{1912}), we have $H^0i_0^{-1}\sK\neq0$ and hence $\sK\neq0$.
\end{proof}

\begin{proof}
For any $1\leq i\leq 3$, $\sQ_i=\sW_i(\C_{T_A^{\rm an}})/\sW_{i-1}(\C_{T_A^{\rm an}})$
\end{proof}>\fi
\begin{exam}\label{hui}
For \[A=\left(\begin{array}[c]{llll}1&0&0&1\\0&1&0&1\\0&0&1&-1\end{array}\right),\]
the canonical epimorphism
$$\alpha_{\C_{T_A^{\rm an}}}\colon \sW_2(\C_{T_A^{\rm an}})/\sW_1(\C_{T_A^{\rm an}})\twoheadrightarrow\bigoplus_{F\in \mathfrak F_2(A)}(\bar j_{F,A})_{!*}(\C_{T_{F}^{\rm an}}[1])$$
is not an isomorphism.
\end{exam}

\begin{proof}
Put $\sQ=\sW_2(\C_{T_A^{\rm an}})/\sW_1(\C_{T_A^{\rm an}})$. Let $i_0\colon {\rm Spec}\,\C\to X_A$ be the closed immersion defined by the unique $T_A$-fixed point $x_0$ on $X_A(\C)$. For any $1\leq k\leq 4$, let $F_k$ be the $k$-th column vector of $A$.
By \cite[8.2.11]{HTT},
\begin{align}\label{s1}&\sW_1(\C_{T_A^{\rm an}})=\tau^{<0}(\bar k_2^A)_*\tau^{<-1}(\bar \ell_2^A)_*(\C_{T_A^{\rm an}}[3]);\\
\label{w6}&\sW_2(\C_{T_A^{\rm an}})=\tau^{<0}(\bar j_A)_*(\C_{T_A^{\rm an}}[3]).\end{align}
This proves that $H^i\sQ=0$ for any $i\geq0$.
By Lemma \ref{1215}, $(\bar k_2^A)^{-1}(\alpha_{\C_{T_A^{\rm an}}})$ is an isomorphism, and the kernel $\sK$ of $\alpha_{\C_{T_A^{\rm an}}}$ is therefore a perverse sheaf on $X_A$ supported on $\{x_0\}=X_A-U_2(A)$. So $H^i\sK=0$ for any $i\neq0$. For any $k$, $\bar j_{F_k}$ is the open immersion $T_{F_k}={\rm Spec}\,\C[x_k^\pm]\hookrightarrow X_{F_k}={\rm Spec}\,\C[x_k]$. Then
$$(\bar j_{F_k,\,A})_{!*}(\C_{T_{F_k}^{\rm an}}[1])=(\bar i_{F_k,\,A})_*(\bar j_{F_k})_{!*}(\C_{T_{F_k}^{\rm an}}[1])=(\bar i_{F_k,\,A})_*\C_{X_{F_k}^{\rm an}}[1].$$
Consequently, $i_0^{-1}(\bar j_{F_k,\,A})_{!*}(\C_{T_{F_k}^{\rm an}}[1])=\C[1]. $ Thus, according to the short exact sequence $$0\to\sK\to\sQ\to\bigoplus_{k=1}^4(\bar j_{F_k,A})_{!*}(\C_{T_{F_k}^{\rm an}}[1])\to0$$ of perverse sheaves, we have $H^i\sQ=0$ for any $i\neq-1$.
Applying $H^i$ to the distinguished triangle $$i_0^{-1}\sK\to i_0^{-1}\sQ\to\bigoplus_{k=1}^4i_0^{-1}(\bar i_{F_k,\,A})_{*}\C_{X_{F_k}}[1],$$ we thus have a short exact sequence of $\C$-vector spaces:
\begin{eqnarray}\label{1911}0\to H^{-1}i_0^{-1}\sQ\to\C^4\to H^0i_0^{-1}\sK\to0.\end{eqnarray}
By (\ref{diao}) and (\ref{w6}), $$H^{-1}i_0^{-1}\sW_2(\C_{T_A^{\rm an}})= H^{-1}i_0^{-1}(\bar j_A)_*(\C_{T_A^{\rm an}}[3])=H^2i_0^{-1}(\bar j_A)_*(\C_{T_A^{\rm an}})=H^2(\textstyle\bigwedge^\bullet\C^3)=\C^3.$$
Applying $H^i$ to the distinguished triangle
$$i_0^{-1}\tau^{<0}(\bar k_2^A)_*\tau^{<-1}(\bar \ell_2^A)_*(\C_{T_A^{\rm an}}[3])\to i_0^{-1}\tau^{<0}(\bar j_A)_*(\C_{T_A^{\rm an}}[3])\to i_0^{-1}\sQ,$$
we obtain an exact sequence
\begin{eqnarray}\label{1912}0\to H^{-1}i_0^{-1}\sW_1(\C_{T_A^{\rm an}})\to H^{-1}i_0^{-1}\sW_2(\C_{T_A^{\rm an}})=\C^3\to H^{-1}i_0^{-1}\sQ\to0.\end{eqnarray}
Combining (\ref{1911}) with (\ref{1912}), we have $H^0i_0^{-1}\sK\neq0$ and hence $\sK\neq0$.
\end{proof}

\bibliographystyle{plain}

\end{document}